\documentclass[11pt,a4paper]{article}

\usepackage{latexsym,amssymb, amsthm}
\usepackage[centertags]{amsmath}
\usepackage{epsfig,pifont}
\usepackage[inactive]{srcltx}  

\usepackage{graphics,helvet,times,mathptm,epic,eepic}

\usepackage{color}

\usepackage{newlfont}

\usepackage{amsfonts,bbm}

\newenvironment{keywords}{ \noindent {\small\bf Key Words}:}{ }

\def\bd{\begin{description}}
\def\ed{\end{description}}

\def\beq{\begin{equation}}
\def\eeq{\end{equation}}
\def\bea{\begin{eqnarray}}
\def\eea{\end{eqnarray}}
\def\beas{\begin{eqnarray*}}
\def\eeas{\end{eqnarray*}}

\newtheorem{lemma}{Lemma}[section]
\newtheorem{theorem}{Theorem}[section]
\newtheorem{corollary}{Corollary}[section]

\theoremstyle{remark}
\newtheorem{example}{Example}[section]

\begin{document}

\title{\textbf{Numerical point of view on Calculus for functions assuming finite, infinite, and infinitesimal values over finite, infinite, and infinitesimal domains}}

\newcommand{\nms}{\normalsize}
\author{\\ {   \bf Yaroslav D. Sergeyev\footnote{Yaroslav D. Sergeyev,
Ph.D., D.Sc., is Distinguished Professor  at the University of
Calabria, Rende, Italy.
 He is also Full Professor (part-time contract) at the N.I.~Lobatchevsky State University,
  Nizhni Novgorod, Russia  and Affiliated Researcher at the Institute of High Performance
Computing and Networking of the National Research Council of
Italy.  }
         }\\ \\ [-2pt]
       \nms   Universit\`a della Calabria,\\[-4pt]
       \nms 87030 Rende (CS)  -- Italy\\ \\[-4pt]
       \nms http://wwwinfo.deis.unical.it/$\sim$yaro\\[-4pt]
         \nms {\tt  yaro@si.deis.unical.it }
}

\date{}

\maketitle

\vspace{-1cm}

\begin{abstract}
The goal of this paper consists  of developing a new (more
physical and numerical in comparison with standard and
non-standard analysis approa\-ches) point of view on Calculus with
functions assuming infinite and infinitesimal values. It uses
recently introduced infinite and infinitesimal numbers being in
accordance with the principle `The part is less than the whole'
observed in the physical world around us. These numbers have a
strong practical advantage with respect to traditional approaches:
they are representable at a new kind of a computer -- the Infinity
Computer -- able to work numerically with all of them. An
introduction to the theory of physical and mathematical continuity
and differentiation (including subdifferentials) for functions
assuming finite, infinite, and infinitesimal values over finite,
infinite, and infinitesimal domains is developed in the paper.
This theory allows one to work with derivatives   that can assume
not only finite but infinite  and infinitesimal values, as well.
It is emphasized that the newly introduced notion of the physical
continuity allows one to see the same mathematical object as a
continuous or a discrete one, in dependence on the wish of the
researcher, i.e., as it happens in the physical world where the
same object can be viewed as a continuous or a discrete in
dependence on the instrument of the observation used by the
researcher. Connections between pure mathematical concepts and
their computational realizations are continuously emphasized
through the text. Numerous examples   are given.
 \end{abstract}



\newpage

\section{Introduction}
\label{s1}

 Numerous trials
have been done during the centuries in order to evolve existing
numeral systems\footnote{ We are reminded that a \textit{numeral}
is a symbol or group of symbols that represents a \textit{number}.
The difference between numerals and numbers is the same as the
difference between words and the things they refer to. A
\textit{number} is a concept that a \textit{numeral} expresses.
The same number can be represented by different numerals. For
example, the symbols `8', `eight', and `VIII' are different
numerals, but they all represent the same number.} in such a way
that infinite and infinitesimal numbers could be included in them
(see \cite{Benci,Cantor,Conway,Leibniz,Newton,Robinson,Wallis}).
Particularly, in the early history of the calculus, arguments
involving infinitesimals played a pivotal role in the derivation
developed by Leibniz and Newton (see \cite{Leibniz,Newton}). The
notion of an infinitesimal, however, lacked a precise mathematical
definition and in order  to provide a more rigorous   foundation
for the calculus infinitesimals were gradually replaced by the
d'Alembert-Cauchy  concept of a limit (see
\cite{Cauchy,DAlembert}).

The creation of a mathematical theory of infinitesimals on which
to base the calculus remained an open problem until the end of the
1950s when Robinson (see \cite{Robinson})   introduced his famous
non-standard Analysis approach. He has shown that non-archimedean
ordered field extensions of the reals contained numbers that could
serve the role of infinitesimals and their reciprocals   could
serve as infinitely large numbers. Robinson then has derived the
theory of limits,  and more generally of Calculus, and has found a
number of important applications of his ideas in many other fields
of Mathematics (see \cite{Robinson}).

In his approach,  Robinson used mathematical tools and terminology
(cardinal numbers, countable sets, continuum, one-to-one
correspondence, etc.) taking their origins from the famous ideas
of Cantor (see \cite{Cantor}) who has shown that there existed
infinite sets having different number of elements. It is well
known nowadays that while dealing with infinite sets, Cantor's
approach leads to some counterintuitive situations that often are
called by non-mathematicians `paradoxes'. For example, the set of
even numbers, $\mathbb{E}$, can be put in a one-to-one
correspondence with the set of all natural numbers, $\mathbb{N}$,
in spite of the fact that $\mathbb{E}$ is  a part of $\mathbb{N}$:
 \beq
\begin{array}{lccccccc}
  \mbox{even numbers:}   & \hspace{5mm} 2, & 4, & 6, & 8,  & 10, & 12, & \ldots    \\

& \hspace{5mm} \updownarrow &  \updownarrow & \updownarrow  & \updownarrow  & \updownarrow  &  \updownarrow &   \\

  \mbox{natural numbers:}& \hspace{5mm}1, &  2, & 3, & 4 & 5,
       & 6,  &    \ldots \\
     \end{array}
\label{4.4.1}
 \eeq
The philosophical principle of Ancient Greeks `\textit{The part is
less than the whole}' observed in the world around us does not
hold true for infinite   numbers   introduced by Cantor, e.g., it
follows $x+1=x$, if $x$ is an infinite cardinal,  although for any
finite $x$ we have   $x+1>x$. As a consequence, the same effects
necessary have reflections in the non-standard Analysis of
Robinson (this is not the case of the interesting non-standard
approach introduced recently in \cite{Benci}).

Due to the enormous importance of the concepts of infinite and
infinitesimal in science, people try to introduce them in their
work with computers, too  (see, e.g. the IEEE Standard for Binary
Floating-Point Arithmetic). However, non-standard Analysis remains
a very theoretical field   because various arithmetics (see
\cite{Benci,Cantor,Conway,Robinson}) developed for infinite and
infinitesimal numbers are quite different with respect to the
finite arithmetic we are used to deal with. Very often they leave
undetermined many operations where infinite numbers take part (for
example, $\infty-\infty$, $\frac{\infty}{\infty}$,  sum of
infinitely many items, etc.) or use representation of infinite
numbers based on infinite sequences of finite numbers. These
crucial difficulties did  not allow people to construct computers
that would be able to work with infinite and infinitesimal numbers
\textit{in the same manner} as we are used to do with finite
numbers  and to study infinite and infinitesimal objects
numerically.

Recently a new applied point of view on infinite and infinitesimal
numbers  has been introduced in \cite{Sergeyev,chaos,informatica}.
The new approach does not use Cantor's ideas and describes
infinite and infinitesimal numbers that are in accordance with the
principle `The part is less than the whole'. It  gives a
possibility to work with finite, infinite, and infinitesimal
quantities \textit{numerically} by using a new kind of a computer
-- the Infinity Computer -- introduced in
\cite{Sergeyev_patent,www,Poland}.  It is worthwhile   noticing
that the new approach   does not contradict Cantor. In contrast,
it can be viewed as an evolution of his deep ideas regarding the
existence of different infinite numbers in a more applied way. For
instance, Cantor has shown that there exist infinite sets having
different cardinalities $\aleph_0$ and $\aleph_1$. In its turn,
the new approach specifies this result showing that in certain
cases within each of these classes it is possible to distinguish
sets with the number of elements being different infinite numbers.

The goal of this paper consists  of developing a new (more
physical and numerical in comparison with standard and
non-standard Analysis approaches) point of view on Calculus. On
the one hand, it uses the approach introduced in
\cite{Sergeyev,chaos,informatica} and, on the other hand, it
incorporates in Calculus the following two main ideas.

i) Note that  foundations of Analysis have been developed more
than 200 years ago with the goal to develop mathematical tools
allowing one to solve problems arising in the real world, as a
result, they reflect  ideas that people had about Physics in that
time. Thus, Analysis that we use now does not include numerous
achievements of Physics of the XX-th century. The brilliant
efforts of Robinson made in the middle of the XX-th century have
been also directed to a reformulation of the classical Analysis in
terms of infinitesimals and not to the creation of a new kind of
Analysis that would incorporate new achievements of Physics. In
fact, he wrote in paragraph 1.1 of his famous book
\cite{Robinson}: `It is shown in this book that Leibniz' ideas can
be fully vindicated and that they lead to a novel and fruitful
approach to classical Analysis and to many other branches of
mathematics'.

The point of view on Calculus   presented in this paper  uses
strongly two methodological ideas borrowed from Physics:
relativity and interrelations holding between the object of an
observation and the tool used for this observation. The latter is
directly related to connections between Analysis and Numerical
Analysis because the numeral systems we use to write down numbers,
functions, etc. are among our tools of investigation and, as a
result, they strongly influence our capabilities to study
mathematical objects.

ii) Both standard and non-standard Analysis mainly study functions
assuming finite values. In this paper, we develop a differential
calculus for functions that can assume finite, infinite, and
infinitesimal values and can be defined over finite, infinite, and
infinitesimal domains. This theory allows one to work with
derivatives   that can assume not only finite but infinite  and
infinitesimal values, as well. Infinite and infinitesimal numbers
are not auxiliary entities in the new Calculus, they are full
members in it and can be used in the same way as finite constants.
In addition, it is important to emphasize that each positive
infinite integer number $a$ expressible in the new numeral system
from \cite{Sergeyev,chaos,informatica} and used in the new
Calculus can be associated with  infinite sets having exactly $a$
elements.

The rest of the paper is structured as follows. In
Section~\ref{s2}, we give a brief introduction to the new
methodology. Section~\ref{s3} describes some preliminary results
dealing with infinite sequences, calculating the number of
elements in various infinite sets, calculating divergent series
and executing arithmetical operations with the obtained infinite
numbers. In Section~\ref{s4}, we introduce two notions of
continuity (working for functions assuming not only finite but
infinite and infinitesimal values, as well) from the points of
view of Physics and Mathematics without usage of the concept of
limit. Section~\ref{s5} describes differential calculus (including
subdifferentials) with functions that can assume finite, infinite,
and infinitesimal values and can be defined over finite, infinite,
and infinitesimal domains. Connections between pure mathematical
concepts and their computational realizations are continuously
emphasized through the text. After all, Section~\ref{s7} concludes
the paper.

We close this Introduction by emphasizing   that the new approach
is introduced as an evolution of standard and non-standard
Analysis and not as a contraposition to them. One or another
version of Analysis can be chosen by the working mathematician in
dependence on the problem he deals with.

\section{Methodology}
\label{s2}

 In this section, we give   a brief
introduction to the new methodology    that can be found in a
rather comprehensive form in the survey \cite{informatica}
downloadable from \cite{www} (see also the monograph
\cite{Sergeyev} written in a popular manner). A number of
applications of the new approach can be found in
\cite{chaos,spirals,biology,Korea}. We start by introducing three
postulates that will fix our methodological positions (having a
strong applied character) with respect to infinite and
infinitesimal quantities and Mathematics, in general.

Usually, when mathematicians deal with infinite objects (sets or
processes) it is supposed   that human beings are able to execute
certain operations infinitely many times (e.g., see
(\ref{4.4.1})). Since we live in a finite world and all human
beings and/or computers finish operations they have started, this
supposition is not adopted.

 \textbf{Postulate 1.} \textit{There  exist
infinite and infinitesimal objects but   human beings and machines
are able to execute only a finite number of operations.}

Due to this Postulate, we accept a priori that we shall never be
able to give a complete description of infinite processes and sets
due to our finite capabilities.

The second postulate is adopted  following the way of reasoning
used in natural sciences where researchers use tools to describe
the object of their study and the instrument used influences the
results of observations. When physicists see a black dot in their
microscope they cannot say: The object of observation \textit{is}
the black dot. They are obliged to say: the lens used in the
microscope allows us to see the black dot and it is not possible
to say anything more about the nature of the object of observation
until we  change the instrument - the lens or the microscope
itself - by a more precise one.

Due to Postulate 1, the same happens in Mathematics studying
natural phenomena, numbers, and objects that can be constructed by
using numbers. Numeral systems used to express numbers are among
the instruments of observations used by mathematicians. Usage of
powerful numeral systems gives the possibility to obtain more
precise results in mathematics in the same way as usage of a good
microscope gives the possibility of obtaining more precise results
in Physics. However, the capabilities of the tools will be always
limited due to Postulate 1 and due to Postulate~2 we shall never
tell, \textbf{what is}, for example, a number but shall just
observe it through numerals expressible in a chosen numeral
system.

 \textbf{Postulate
2.} \textit{We shall not   tell \textbf{what are} the mathematical
objects we deal with; we just shall construct more powerful tools
that will allow us to improve our capacities to observe and to
describe properties of mathematical objects.}

Particularly, this means that from our point of view, axiomatic
systems do not define mathematical objects but just determine
formal rules for operating with certain numerals reflecting some
properties of the studied mathematical objects. Throughout the
paper, we shall always emphasize this philosophical triad --
researcher, object of investigation, and tools used to observe the
object -- in various mathematical and computational contexts.

Finally, we adopt the principle of Ancient Greeks mentioned above
as  the third postulate.

\textbf{Postulate 3.} \textit{The principle `The part is less than
the whole' is applied to all numbers (finite, infinite, and
infinitesimal) and to all sets and processes (finite and
infinite).}

Due to this declared applied statement, it becomes clear that the
subject of this paper is out of Cantor's approach and, as a
consequence, out of non-standard analysis of Robinson. Such
concepts as bijection, numerable and continuum sets, cardinal and
ordinal numbers cannot be used in this paper because they belong
to the theory working with different assumptions. However, the
approach used here does not contradict Cantor and Robinson. It can
be viewed just as a more strong lens of a mathematical microscope
that allows one to distinguish more objects and to work with them.

In \cite{Sergeyev,informatica}, a
  new  numeral system has
been developed in accordance with Postulates 1--3. It gives one a
possibility to execute numerical computations not only with finite
numbers but also with infinite and infinitesimal ones. The main
idea consists of the possibility to measure infinite and
infinitesimal quantities by different (infinite, finite, and
infinitesimal) units of measure.

A new infinite unit of measure   has been introduced for this
purpose   as the number of elements of the set $\mathbb{N}$ of
natural numbers. It is expressed by the numeral \ding{172} called
\textit{grossone}. It is necessary to note immediately that
\ding{172} is neither Cantor's $\aleph_0$ nor $\omega$.
Particularly, it has both cardinal and ordinal properties as usual
finite natural numbers (see \cite{informatica}).

Formally, grossone is introduced as a new number by describing its
properties postulated by the \textit{Infinite Unit Axiom} (IUA)
(see \cite{Sergeyev,informatica}). This axiom is added to axioms
for real numbers similarly to addition of the axiom determining
zero to axioms of natural numbers when integer numbers are
introduced. It is important to emphasize that we speak about
axioms of real numbers in sense of Postulate~2, i.e., axioms
define formal rules of operations with numerals in a given numeral
system.

Inasmuch as it has been postulated that grossone is a number,  all
other axioms for numbers hold for it, too. Particularly,
associative and commutative properties of multiplication and
addition, distributive property of multiplication over addition,
existence of   inverse  elements with respect to addition and
multiplication hold for grossone as for finite numbers. This means
that  the following relations hold for grossone, as for any other
number
 \beq
 0 \cdot \mbox{\ding{172}} =
\mbox{\ding{172}} \cdot 0 = 0, \hspace{3mm}
\mbox{\ding{172}}-\mbox{\ding{172}}= 0,\hspace{3mm}
\frac{\mbox{\ding{172}}}{\mbox{\ding{172}}}=1, \hspace{3mm}
\mbox{\ding{172}}^0=1, \hspace{3mm}
1^{\mbox{\tiny{\ding{172}}}}=1, \hspace{3mm}
0^{\mbox{\tiny{\ding{172}}}}=0.
 \label{3.2.1}
       \eeq

Let us comment upon the nature of grossone by some illustrative
examples.

\begin{example}
\label{e1}  Infinite numbers constructed using grossone  can be
interpreted in terms of the number of elements of infinite sets.
For example, $\mbox{\ding{172}}-1$ is the number of elements of a
set $B=\mathbb{N}\backslash\{b\}$, $b \in \mathbb{N}$, and
$\mbox{\ding{172}}+1$ is the number of elements of a set
$A=\mathbb{N}\cup\{a\}$, where $a \notin \mathbb{N}$.  Due to
Postulate~3, integer positive numbers that are larger than
grossone do not belong to $\mathbb{N}$ but also can be easily
interpreted. For instance,   $\mbox{\ding{172}}^2$ is the number
of elements of the set
  $V$, where  $ V  =
\{ (a_1, a_2)  : a_1 \in   \mathbb{N}, a_2 \in   \mathbb{N}  \}.
 $   \hfill $\Box$

\end{example}

\begin{example}
\label{e2}  Grossone has been introduced as the quantity of
natural numbers. As a consequence, similarly  to the set
 \beq
  A=\{1, 2, 3, 4, 5\}
\label{4.1.deriva_0}
 \eeq
   consisting of
5 natural numbers where 5 is the largest number in $A$, \ding{172}
is the largest    number\footnote{This fact is one of the
important methodological differences with respect to non-standard
analysis theories where it is supposed that infinite numbers   do
not belong to $\mathbb{N}$.} in $\mathbb{N}$ and
$\mbox{\ding{172}} \in \mathbb{N}$ analogously to the fact that 5
belongs to $A$. Thus, the set, $\mathbb{N}$, of natural numbers
can be written  in the form
 \beq
\mathbb{N} = \{ 1,2,  \hspace{3mm} \ldots  \hspace{3mm}
\frac{\mbox{\ding{172}}}{2}-2, \frac{\mbox{\ding{172}}}{2}-1,
\frac{\mbox{\ding{172}}}{2}, \frac{\mbox{\ding{172}}}{2}+1,
\frac{\mbox{\ding{172}}}{2}+2, \hspace{3mm}  \ldots \hspace{3mm}
\mbox{\ding{172}}-2, \hspace{2mm}\mbox{\ding{172}}-1, \hspace{2mm}
\mbox{\ding{172}} \}.   \label{4.1}
       \eeq
Note that traditional numeral systems did not allow us to see
infinite natural numbers
 \beq \ldots  \hspace{3mm}
\frac{\mbox{\ding{172}}}{2}-2, \frac{\mbox{\ding{172}}}{2}-1,
\frac{\mbox{\ding{172}}}{2}, \frac{\mbox{\ding{172}}}{2}+1,
\frac{\mbox{\ding{172}}}{2}+2, \hspace{3mm} \ldots  \hspace{3mm}
\mbox{\ding{172}}-2, \mbox{\ding{172}}-1, \mbox{\ding{172}}.
\label{4.1.deriva_1}
 \eeq
Similarly,
  Pirah\~{a}\footnote{Pirah\~{a} is a primitive tribe living in
Amazonia  that uses a very simple numeral system for counting:
one, two, `many'(see \cite{Gordon}). For Pirah\~{a}, all
quantities larger than two are just `many' and such operations as
2+2 and 2+1 give the same result, i.e., `many'. Using their weak
numeral system Pirah\~{a} are not able to distinguish numbers
larger than 2 and, as a result, to execute arithmetical operations
with them. Another peculiarity of this numeral system  is that
\mbox{`many'}+ 1= \mbox{`many'}. It can be immediately seen that
this result  is very similar to our traditional record $\infty +
1= \infty$.}     are not able to see  finite numbers larger than 2
using their weak numeral system but these numbers are visible if
one uses a more powerful numeral system. Due to Postulate~2, the
same object  of observation -- the set $\mathbb{N}$ --   can be
observed by different instruments -- numeral systems -- with
different accuracies allowing one to express  more or less natural
numbers. \hfill $\Box$
\end{example}

This example illustrates also the fact that when we speak about
sets (finite or infinite) it is necessary to take care about tools
used to describe a set (remember Postulate~2). In order to
introduce a set, it is necessary to have a language (e.g., a
numeral system) allowing us to describe its elements and the
number of the elements in the set. For instance, the set $A$ from
(\ref{4.1.deriva_0}) cannot be defined using the mathematical
language of Pirah\~{a}.

Analogously, the words `the set of all finite numbers' do not
define a set completely from our point of view, as well. It is
always necessary to specify which instruments are used to describe
(and to observe) the required set and, as a consequence, to speak
about `the set of all finite numbers expressible in a fixed
numeral system'. For instance, for Pirah\~{a} `the set of all
finite numbers'  is the set $\{1, 2 \}$ and for another Amazonian
tribe -- Munduruk\'u\footnote{Munduruk\'u (see \cite{Pica}) fail
in exact arithmetic with numbers larger than   5 but are able to
compare and add large approximate numbers that are far beyond
their naming range. Particularly, they use the words `some, not
many' and `many, really many' to distinguish two types of large
numbers (in this connection think about Cantor's  $\aleph_0$ and
$\aleph_1$).} -- `the set of all finite numbers' is the set $A$
from (\ref{4.1.deriva_0}). As it happens in Physics, the
instrument used for an observation bounds the possibility of the
observation. It is not possible to say how we shall see the object
of our observation if we have not clarified which instruments will
be used to execute the observation.

Introduction of  grossone gives us a possibility to compose new
(in comparison with traditional numeral systems) numerals and to
see through them not only numbers (\ref{4.1.deriva_0}) but also
certain numbers larger than \ding{172}. We can speak about the set
of \textit{extended natural numbers} (including $\mathbb{N}$ as a
proper subset) indicated as $\widehat{\mathbb{N}}$ where
 \beq
  \widehat{\mathbb{N}} = \{
1,2, \ldots ,\mbox{\ding{172}}-1, \mbox{\ding{172}},
\mbox{\ding{172}}+1, \mbox{\ding{172}}+2, \mbox{\ding{172}}+3,
\ldots , \mbox{\ding{172}}^2-1, \mbox{\ding{172}}^2.
\mbox{\ding{172}}^2+1, \ldots \} \label{4.2.2}
       \eeq
However, analogously to the situation with `the set of all finite
numbers', the number of elements of the set $\widehat{\mathbb{N}}$
cannot be expressed within a numeral system using only \ding{172}.
It is necessary to introduce in a reasonable way a more powerful
numeral system and to define   new numerals (for instance,
\ding{173}, \ding{174}, etc.) of this system that would allow one
to fix the set (or sets) somehow. In general, due to Postulate~1
and~2, for any fixed numeral $\mathcal{A}$ system there always be
sets that cannot be described using $\mathcal{A}$.

\begin{example}
\label{e3}

Analogously to (\ref{4.1}), the set, $\mathbb{E}$, of even natural
numbers can be written now in the form
 \beq
\mathbb{E} = \{ 2,4,6 \hspace{5mm} \ldots  \hspace{5mm}
\mbox{\ding{172}}-4, \hspace{2mm}\mbox{\ding{172}}-2, \hspace{2mm}
\mbox{\ding{172}} \}.   \label{4.1.0}
       \eeq
Due to Postulate 3 and the IUA (see \cite{Sergeyev,informatica}),
it follows that the number of elements of the set of even numbers
is equal to $\frac{\mbox{\ding{172}}}{2}$ and \ding{172} is even.
Note that the next even number is $\mbox{\ding{172}}+2$ but it is
not natural because $\mbox{\ding{172}}+2  > \mbox{\ding{172}}$, it
is extended natural (see (\ref{4.2.2})). Thus, we can write down
not only initial (as it is  done traditionally) but also the final
part of (\ref{4.4.1})
  \[
\begin{array}{cccccccccc}
 2, & 4, & 6, & 8,  & 10, & 12, & \ldots  &
\mbox{\ding{172}} -4,  &    \mbox{\ding{172}}  -2,   &    \mbox{\ding{172}}    \\
 \updownarrow &  \updownarrow & \updownarrow  &
\updownarrow  & \updownarrow  &  \updownarrow  & &
  \updownarrow    & \updownarrow   &
  \updownarrow
   \\
 1, &  2, & 3, & 4 & 5, & 6,   &   \ldots  &    \frac{\mbox{\ding{172}}}{2} - 2,  &
     \frac{\mbox{\ding{172}}}{2} - 1,  &    \frac{\mbox{\ding{172}}}{2}   \\
     \end{array}
\]
concluding so (\ref{4.4.1})   in a complete accordance with
Postulate~3. It is worth  noticing that the new numeral system
allows us to solve many other `paradoxes' related to infinite and
infinitesimal quantities (see \cite{Sergeyev,informatica,Korea}).
 \hfill $\Box$
\end{example}

In order to  express numbers having finite, infinite, and
infinitesimal parts, records similar to traditional positional
numeral systems can be used (see \cite{Sergeyev,informatica}). To
construct a number $C$ in the new numeral positional system with
base \ding{172}, we subdivide $C$ into groups corresponding to
powers of \ding{172}:
 \beq
  C = c_{p_{m}}
\mbox{\ding{172}}^{p_{m}} +  \ldots + c_{p_{1}}
\mbox{\ding{172}}^{p_{1}} +c_{p_{0}} \mbox{\ding{172}}^{p_{0}} +
c_{p_{-1}} \mbox{\ding{172}}^{p_{-1}}   + \ldots   + c_{p_{-k}}
 \mbox{\ding{172}}^{p_{-k}}.
\label{3.12}
       \eeq
 Then, the record
 \beq
  C = c_{p_{m}}
\mbox{\ding{172}}^{p_{m}}    \ldots   c_{p_{1}}
\mbox{\ding{172}}^{p_{1}} c_{p_{0}} \mbox{\ding{172}}^{p_{0}}
c_{p_{-1}} \mbox{\ding{172}}^{p_{-1}}     \ldots c_{p_{-k}}
 \mbox{\ding{172}}^{p_{-k}}
 \label{3.13}
       \eeq
represents  the number $C$, where all numerals $c_i\neq0$, they
belong to a traditional numeral system and are called
\textit{grossdigits}. They express finite positive or negative
numbers and show how many corresponding units
$\mbox{\ding{172}}^{p_{i}}$ should be added or subtracted in order
to form the number $C$.

Numbers $p_i$ in (\ref{3.13}) are  sorted in the decreasing order
with $ p_0=0$
\[
p_{m} >  p_{m-1}  > \ldots    > p_{1} > p_0 > p_{-1}  > \ldots
p_{-(k-1)}  >   p_{-k}.
 \]
They are called \textit{grosspowers} and they themselves can be
written in the form (\ref{3.13}).
 In the record (\ref{3.13}), we write
$\mbox{\ding{172}}^{p_{i}}$ explicitly because in the new numeral
positional system  the number   $i$ in general is not equal to the
grosspower $p_{i}$. This gives the possibility to write down
numerals without indicating grossdigits equal to zero.

The term having $p_0=0$ represents the finite part of $C$ because,
due to (\ref{3.2.1}), we have $c_0 \mbox{\ding{172}}^0=c_0$. The
terms having finite positive gross\-powers represent the simplest
infinite parts of $C$. Analogously, terms   having   negative
finite grosspowers represent the simplest infinitesimal parts of
$C$. For instance, the  number
$\mbox{\ding{172}}^{-1}=\frac{1}{\mbox{\ding{172}}}$ is
infinitesimal. It is the inverse element with respect to
multiplication for \ding{172}:
 \beq
\mbox{\ding{172}}^{-1}\cdot\mbox{\ding{172}}=\mbox{\ding{172}}\cdot\mbox{\ding{172}}^{-1}=1.
 \label{3.15.1}
       \eeq
Note that all infinitesimals are not equal to zero. Particularly,
$\frac{1}{\mbox{\ding{172}}}>0$ because it is a result of division
of two positive numbers. All of the numbers introduced above can
be grosspowers, as well, giving thus a possibility to have various
combinations of quantities and to construct  terms having a more
complex structure.

\begin{example}
\label{e4} In this example, it is shown   how to write down
numbers in the new numeral system and   how the value of the
number is calculated:
\[
 C_1=7.6\mbox{\ding{172}}^{24.5\mbox{\small{\ding{172}}}-7.1}\,34\mbox{\ding{172}}^{3.2}({\mbox{\small{-}}}3)\mbox{\ding{172}}^{\mbox{\tiny\ding{172}}^{-1}}70\mbox{\ding{172}}^{0}52.1\mbox{\ding{172}}^{-6.8}({\mbox{\small{-}}}0.23)\mbox{\ding{172}}^{-9.4\mbox{\tiny\ding{172}}}
=
 \]
 \[
7.6\mbox{\ding{172}}^{24.5\mbox{\small{\ding{172}}}-7.1}\,+34\mbox{\ding{172}}^{3.2}-3\mbox{\ding{172}}^{\mbox{\tiny\ding{172}}^{-1}}+70\mbox{\ding{172}}^{0}+52.1\mbox{\ding{172}}^{-6.8}-0.23\mbox{\ding{172}}^{-9.4\mbox{\tiny\ding{172}}}.
\]
The number $C_1$ above has two infinite parts of the type
$\mbox{\ding{172}}^{24.5\mbox{\small{\ding{172}}}-7.1}$ and
$\mbox{\ding{172}}^{3.2}$, one part
$\mbox{\ding{172}}^{\mbox{\tiny\ding{172}}^{-1}}$ that is
infinitesimally close to $\mbox{\ding{172}}^{0}$, a finite part
corresponding to $\mbox{\ding{172}}^{0}$, and two infinitesimal
parts of the type $\mbox{\ding{172}}^{-6.8}$ and
$\mbox{\ding{172}}^{-9.4\mbox{\tiny\ding{172}}}$. The
corresponding grossdigits show how many units of each kind should
be taken (added or subtracted) to form~$C_1$. \hfill$\Box$
\end{example}

\section{Preliminary results} \label{s3}

\subsection{Infinite sequences}
\label{s3.1}

We start by recalling   traditional definitions of the infinite
sequences and subsequences.  An \textit{infinite sequence}
$\{a_n\}, a_n \in A, n \in \mathbb{N},$ is a function having as
the domain the set of natural numbers, $\mathbb{N}$, and as the
codomain  a set $A$. A \textit{subsequence} is   a sequence from
which some of its elements have been removed.

Let us look at these definitions from the new point of view.
Grossone has been introduced as the number of elements of the set
$\mathbb{N}$. Thus, due to the sequence definition given above,
any sequence having $\mathbb{N}$ as the domain  has \ding{172}
elements.

The notion of subsequence is introduced as a sequence from which
some of its elements have been removed. The new numeral system
gives the possibility to indicate explicitly the removed elements
and to count how many they are  independently of the fact whether
their numbers in the sequence are finite or infinite. Thus, this
  gives infinite sequences having a number of members
less than~\ding{172}. Then the following result holds.
\begin{theorem}
\label{t2} The number of elements of any infinite sequence is less
or equal to~\ding{172}.
\end{theorem}

\textit{Proof.}  The proof is obvious and is so omitted.  \hfill
$\Box$

One of the immediate consequences of the understanding of this
result is that any sequential process can have at maximum
\ding{172} elements and, due to Postulate 1, it depends on the
chosen numeral system which numbers among  \ding{172} members of
the process we can observe (see \cite{Sergeyev,chaos} for a
detailed discussion). Particularly, this means that from a set
having more than grossone elements it is not possible to choose
all its elements if only one sequential process of choice is used
for this purpose. Another important thing that we can do now with
the infinite sequence is the possibility to observe their final
elements if they are expressible in the chosen numeral system, in
the same way as it happens with finite sequences.

It becomes appropriate now to define the \textit{complete
sequence} as an infinite sequence  containing \ding{172} elements.
For example, the sequence $\{1, 2, 3, \ldots \mbox{\ding{172}}-2,
\mbox{\ding{172}}-1, \mbox{\ding{172}} \}$ of natural numbers is
complete, the sequences  $\{2, 4, 6, \ldots \mbox{\ding{172}}-4,
\mbox{\ding{172}}-2, \mbox{\ding{172}} \}$ and $\{1, 3, 5, \ldots
\mbox{\ding{172}}-5, \mbox{\ding{172}}-3, \mbox{\ding{172}}-1 \}$
of even  and odd natural numbers  are not complete because, due to
the IUA (see \cite{Sergeyev,chaos}), they have
$\frac{\mbox{\ding{172}}}{2}$ elements each. Thus, to describe a
sequence we should use the record $\{a_n: k \}$ where $a_n$ is, as
usual, the general element and $k$ is the number (finite or
infinite) of members of the sequence.

\begin{example}
\label{e25} Let us consider the set, $\widehat{\mathbb{N}}$, of
extended natural numbers from (\ref{4.2.2}). Then, starting from
the number 3, the process of the sequential counting can arrive at
maximum to $\mbox{\ding{172}}+2$:
\[
 1,2,\underbrace{3,4,\hspace{1mm}  \ldots \hspace{1mm}
\mbox{\ding{172}}-2,\hspace{1mm}
 \mbox{\ding{172}}-1,
\mbox{\ding{172}},  \mbox{\ding{172}}+1,
\mbox{\ding{172}}+2}_{\mbox{\ding{172}}}, \mbox{\ding{172}}+3,
  \ldots
\]
Analogously, starting from the number
$\frac{\mbox{\ding{172}}}{2}+1$, the following process of the
sequential counting
\[
 \ldots  \hspace{1mm}
  \frac{\mbox{\ding{172}}}{2}-1,
\frac{\mbox{\ding{172}}}{2},
\underbrace{\frac{\mbox{\ding{172}}}{2}+1,
\frac{\mbox{\ding{172}}}{2}+2, \hspace{1mm} \ldots  \hspace{1mm}
\mbox{\ding{172}}-1, \mbox{\ding{172}}, \mbox{\ding{172}}+1,
\hspace{1mm} \ldots \frac{3\mbox{\ding{172}}}{2}-1,
\frac{3\mbox{\ding{172}}}{2}}_{\mbox{\ding{172}}  },
\frac{3\mbox{\ding{172}}}{2}+1,
 \hspace{1mm}  \ldots
\]
can arrive as a maximum to the number
$\frac{3\mbox{\ding{172}}}{2}$. \hfill $\Box$
 \end{example}

\subsection{Series}
\label{s3.2}

Postulate~3 imposes us the same behavior in relation to finite and
infinite quantities. Thus, working with sums it is always
necessary to indicate explicitly the number of items (finite or
infinite) in the sum. Of course, to calculate a sum numerically it
is necessary that the number of items and the result are
expressible in the numeral system used for calculations. It is
important to notice that even though a sequence cannot have more
than \ding{172} elements, the number of items in a sum can be
greater than grossone because the process of summing up  should
not necessarily be executed by a sequential adding of items.

\begin{example}
\label{e6} Let us consider two infinite series $S_1=7+7+7+\ldots$
and $S_2=3+3+3+\ldots$  Traditional analysis gives us a very poor
answer that both of them diverge to infinity. Such operations as
$S_2 - S_1$ or $ \frac{S_1}{S_2} $ are not defined. In our
approach, it is necessary to indicate explicitly the number of
items in the sum and it is not important whether it is  finite or
infinite.

Suppose that the   series $S_1$ has $k$ items and $S_2$ has $n$
items:
$$S_1(k)=\underbrace{7+7+7+\ldots+7}_k, \hspace{1cm} S_2(n)=\underbrace{3+3+3+\ldots+3}_n.$$
Then $S_1(k)=7k$ and $S_2(n)=3n$ and by giving different numerical
values (finite or infinite) to $k$ and $n$ we obtain different
numerical values for the sums.  For chosen $k$ and $n$ it becomes
possible to calculate $S_2(n) - S_1(k)$ (analogously, the
expression $ \frac{S_1(k)}{S_2(n)} $ can be calculated). If, for
instance, $k=5\mbox{\ding{172}}$ and $n=\mbox{\ding{172}}$   we
obtain $S_1(5\mbox{\ding{172}})=35\mbox{\ding{172}}$,
$S_2(\mbox{\ding{172}})=3\mbox{\ding{172}}$ and it follows
\[
S_2(\mbox{\ding{172}}) - S_1(5\mbox{\ding{172}})=
3\mbox{\ding{172}} - 35\mbox{\ding{172}} = - 32\mbox{\ding{172}}<
0.
\]
If  $k=3\mbox{\ding{172}}$ and $n=7\mbox{\ding{172}}+2$   we
obtain $S_1(3\mbox{\ding{172}})=21\mbox{\ding{172}}$,
$S_2(7\mbox{\ding{172}}+2)=21\mbox{\ding{172}}+6$ and  it follows
\[
S_2(7\mbox{\ding{172}}+2) -  S_1(3\mbox{\ding{172}})=
21\mbox{\ding{172}}+6 - 21\mbox{\ding{172}} = 6.
\]
 It is
also possible to sum up sums having an infinite number of infinite
or infinitesimal items
\[
S_3(l)=\underbrace{2\mbox{\ding{172}}+2\mbox{\ding{172}}+\ldots+2\mbox{\ding{172}}}_l,
\hspace{1cm}
S_4(m)=\underbrace{4\mbox{\ding{172}}^{-1}+4\mbox{\ding{172}}^{-1}+\ldots+4\mbox{\ding{172}}^{-1}}_m.
\]
For $l=m=0.5\mbox{\ding{172}}$ it follows
$S_3(0.5\mbox{\ding{172}})=  \mbox{\ding{172}}^2$ and
$S_4(0.5\mbox{\ding{172}})=2$ (recall that
$\mbox{\ding{172}}\cdot\mbox{\ding{172}}^{-1}=\mbox{\ding{172}}^{0}=1$
(see (\ref{3.15.1})).  It can be seen from this example that it is
possible to obtain finite   numbers as the result of summing up
infinitesimals. This is a direct consequence of Postulate~3.
 \hfill
 $\Box$
 \end{example}

The infinite and infinitesimal numbers allow us to  also calculate
arithmetic and geometric series with an infinite number of items.
Traditional approaches tell us that if  $a_{n} = a_{1} + (n - 1)d$
then for a finite $n$ it is possible to use the formula
\[
\sum_{i=1}^{n} a_{i} = \frac{n}{2}(a_{1} + a_{n}).
\]
Due to Postulate~3, we can use  it also for infinite $n$.

\begin{example}
\label{e7} The sum of all natural numbers from 1 to \ding{172} can
be calculated as follows
 \beq
1+2+3+ \ldots + (\mbox{\ding{172}}-1) + \mbox{\ding{172}} =
\sum_{i=1}^{\mbox{\ding{172}}} i = \frac{\mbox{\ding{172}}}{2}(1 +
\mbox{\ding{172}})= 0.5\mbox{\ding{172}}^{2}0.5\mbox{\ding{172}}.
 \label{der1}
 \eeq
 Let us  now calculate the following sum of
infinitesimals where each item is \ding{172} times less than the
corresponding item of (\ref{der1})
\[
\mbox{\ding{172}}^{-1}+2\mbox{\ding{172}}^{-1}+ \ldots +
(\mbox{\ding{172}}-1)\cdot\mbox{\ding{172}}^{-1} +
\mbox{\ding{172}}\cdot\mbox{\ding{172}}^{-1} =
\sum_{i=1}^{\mbox{\ding{172}}} i\mbox{\ding{172}}^{-1} =
\frac{\mbox{\ding{172}}}{2}(\mbox{\ding{172}}^{-1} + 1)=
0.5\mbox{\ding{172}}^{1}0.5.
\]
Obviously, the obtained number, $0.5\mbox{\ding{172}}^{1}0.5$  is
\ding{172} times less than the sum in (\ref{der1}).  This example
shows, particularly, that infinite numbers can also be obtained as
the result of summing up infinitesimals.
 \hfill
 $\Box$
 \end{example}
Let us now consider the geometric series $\sum_{i=0}^{\infty}
q^i$. Traditional analysis proves that it converges to
$\frac{1}{1-q}$ for $q$ such that $-1 < q < 1$. We are able to
give a more precise answer for \textit{all} values of $q$. To do
this we should fix the number of items in the sum. If we suppose
that it contains $n$ items, then
 \beq
 Q_n = \sum_{i=0}^{n}
q^i = 1 + q + q^2 + \ldots + q^n.
 \label{3.7.2.f}
 \eeq
By multiplying the left-hand and the right-hand parts of this
equality by $q$ and by subtracting the result from (\ref{3.7.2.f})
we obtain
\[
Q_n - qQ_n = 1-q^{n+1}
\]
and, as a consequence, for all $q\neq 1$ the formula
 \beq
  Q_n =
(1-q^{n+1})(1-q)^{-1}
 \label{3.7.2.f.1}
 \eeq
holds for finite and infinite $n$. Thus, the possibility to
express infinite and infinitesimal numbers allows us   to take
into account    infinite $n$   and the value $q^{n+1}$ being
infinitesimal for    $q, \,\,-1<q<1$. Moreover, we can calculate
$Q_n$ for infinite and finite values of $n$ and $q=1$, because in
this case we have just
\[
Q_n = \underbrace{1+1+1+\ldots+1}_{n+1} = n+1.
\]
\begin{example}
\label{e8}
 As the first example we consider the divergent series
\[
1 + 3 + 9 + \ldots = \sum_{i=0}^{\infty} 3^i.
\]
To fix it, we should decide the number of items, $n$, at the sum
and, for example, for $n=\mbox{\ding{172}}^2$ we obtain
\[
\sum_{i=0}^{\mbox{\ding{172}}^2} 3^i = 1 + 3 + 9 + \ldots +
3^{\mbox{\ding{172}}^2}= \frac{1-3^{\mbox{\ding{172}}^2+1}}{1-3} =
0.5(3^{\mbox{\ding{172}}^2+1}-1).
\]
Analogously, for  $n=\mbox{\ding{172}}^2+1$ we obtain
\[
1 + 3 + 9 + \ldots + 3^{\mbox{\ding{172}}^2} +
3^{\mbox{\ding{172}}^{2}+1}= 0.5(3^{\mbox{\ding{172}}^2+2}-1).
\]
If we now find the difference between the two sums
\[
0.5(3^{\mbox{\ding{172}}^2+2}-1)
-(0.5(3^{\mbox{\ding{172}}^2+1}-1)) = 3^{\mbox{\ding{172}}^2+1}
(0.5\cdot3 - 0.5) = 3^{\mbox{\ding{172}}^{2}+1}
\]
 we  obtain the
newly added item $3^{\mbox{\ding{172}}^2+1}$. \hfill
 $\Box$
 \end{example}

\begin{example}
\label{e9} In this example, we consider the series
$\sum_{i=1}^{\infty}\frac{1}{2^i}$. It is well-known that it
converges to one. However, we are able to give a more precise
answer. In fact, due to Postulate~3, the formula
 \[
 \sum_{i=1}^{n}\frac{1}{2^i}= \frac{1}{2}( 1+ \frac{1}{2} + \frac{1}{2^2} + \ldots + \frac{1}{2^{n-1}}) = \frac{1}{2} \cdot \frac{1-\frac{1}{2^{n}}}{1-\frac{1}{2}}  = 1-\frac{1}{2^n}
 \]
can be used directly for   infinite $n$, too. For example, if
$n=\mbox{\ding{172}}$ then
$$\sum_{i=1}^{\mbox{\small{\ding{172}}}}\frac{1}{2^i}=1-\frac{1}{2^{\mbox{\tiny{\ding{172}}}}},$$
where $\frac{1}{2^{\mbox{\tiny{\ding{172}}}}}$ is infinitesimal.
Thus, the traditional answer $\sum_{i=1}^{\infty}\frac{1}{2^i}=1$
was just a finite approximation to our more precise result using
infinitesimals. \hfill
 $\Box$
 \end{example}

\subsection{From limits to  expressions}
\label{s3.3}

 Let us now discuss the  theory of limits from the point of
view of our approach. In traditional analysis, if a limit $\lim_{x
\rightarrow a}f(x)$ exists, then it gives us a very poor -- just
one  value -- information about the behavior of $f(x)$ when $x$
tends to $a$. Now we can obtain significantly more rich
information because we are able to calculate $f(x)$ directly at
any finite, infinite, or infinitesimal point that can be expressed
by the new positional system. This can be done even in the cases
where the limit does not exist. Moreover, we can easily work with
functions assuming infinite  or infinitesimal values at infinite
or infinitesimal points.

Thus, limits $\lim_{x \rightarrow \infty}f(x)$ equal to infinity
can be substituted by  precise infinite numerals that are
different for different infinite values of $x$. If we speak about
limits of sequences, $\lim_{n \rightarrow \infty}a(n)$, then $ n
\in \mathbb{N}$ and, as a consequence, it follows from
Theorem~\ref{t2} that $n$ at which we can evaluate $a(n)$ should
be less than or equal to grossone.

\begin{example}
\label{e10} From the traditional point of view, the following two
limits
 \[
\lim_{x \rightarrow +\infty}(7x^8+2x^3)= +\infty, \hspace{1cm}
\lim_{x \rightarrow +\infty}(7x^8+2x^3+10^{100})=
 +\infty.
 \]
give  us   the same result, $+\infty$,  in spite of the fact that
for any finite $x$ it follows
 \[
7x^8+2x^3+10^{100} - (7x^8+2x^3) = 10^{100}
\]
that is a rather huge number. In other words, the two expressions
that are comparable for any finite $x$ cannot be compared at
infinity. The new approach allows us to calculate exact values of
both expressions, $7x^8+2x^3$ and $7x^8+2x^3+10^{100}$, at any
infinite $x$ expressible in the chosen numeral system. For
instance, the choice $x=\mbox{\ding{172}}^{2}$ gives the value
\[
7(\mbox{\ding{172}}^{2})^{8}+2(\mbox{\ding{172}}^{2})^{3}=
7\mbox{\ding{172}}^{16}2\mbox{\ding{172}}^{6}
\]
for the first expression and
$7\mbox{\ding{172}}^{16}2\mbox{\ding{172}}^{6}10^{100}$ for the
second one. We can easily calculate their difference  that
evidently is equal to $10^{100}$.  \hfill
 $\Box$
 \end{example}

Limits with  the argument tending to zero can be considered
analogously. In this case, we can calculate the corresponding
expression at   infinitesimal points using the new positional
system and to obtain   significantly more reach information. If
the traditional limit exists, it will be just a finite
approximation of the new more precise result   having the finite
part and eventual infinitesimal parts.

\begin{example}
\label{e11} Let us consider the following limit
 \beq
 \lim_{h \rightarrow 0}\frac{(3+h)^2-3^2}{h}= 6.
\label{4.5}
 \eeq
In the new positional system for $h\neq0$ we obtain
 \beq
\frac{(3+h)^2-3^2}{h}= 6 + h.
 \label{4.6}
 \eeq
If, for instance, the number $h=\mbox{\ding{172}}^{-1}$, the
answer is $6\mbox{\ding{172}}^{0}\mbox{\ding{172}}^{-1}$, if
$h=4\mbox{\ding{172}}^{-2}$ we obtain
$6\mbox{\ding{172}}^{0}4\mbox{\ding{172}}^{-2}$, etc. Thus, the
value of  the limit (\ref{4.5})  is just the finite approximation
of the number (\ref{4.6}) having finite and infinitesimal parts
that can be used in possible further calculations if an accuracy
higher than the finite one is required. \hfill $\Box$
\end{example}

The new numeral system allows us to  evaluate expressions at
infinite   or infinitesimal points   when their limits do not
exist giving thus a very powerful tool for studying divergent
processes. Another important feature of the new approach consists
of the possibility to construct expressions where infinite and/or
infinitesimal quantities are involved and to evaluate them at
infinite or infinitesimal points.

\begin{example}
\label{e12} Let us consider the following expression
$\frac{1}{\mbox{\ding{172}}}x^2+\mbox{\ding{172}}x +2$. For
example, for the infinite $x=3\mbox{\ding{172}}$ we obtain the
infinite value $9\mbox{\ding{172}} +
3\mbox{\ding{172}}^{2}+2=3\mbox{\ding{172}}^{2}9\mbox{\ding{172}}^{1}2$.
For the infinitesimal $x=\mbox{\ding{172}}^{-1}$ we have
$\mbox{\ding{172}}^{-3} +
1+2=3\mbox{\ding{172}}^{0}\mbox{\ding{172}}^{-3}$.
 \hfill
$\Box$
\end{example}

\subsection{Expressing and counting points over one-dimensional intervals}
\label{s3.4}

We start this subsection    by calculating the number of points at
the interval $[0,1)$. To do this we need a definition of the term
`point'\index{point} and mathematical tools to indicate a point.
Since this concept is one of the most fundamental, it is very
difficult to find an adequate definition for it. If we accept (as
is usually done in modern mathematics) that  a \textit{point}  in
$[0,1)$ is determined by a numeral $x$   called the
\textit{coordinate of the point} where $x \in \mathcal{S}$ and
 $\mathcal{S}$ is a set of numerals,    then we can
indicate  the point  by its coordinate  $x$  and are able to
execute required calculations.

It is important to emphasize  that we have not postulated that $x$
belongs to the  set, $\mathbb{R}$, of real numbers  as it is
usually done.  Since we can express coordinates only by numerals,
then different choices of numeral systems lead to various sets of
numerals and, as a consequence, to different sets of points we can
refer to.  The choice of a numeral system will define what is the
\textit{point} for us and we shall not be able to work with those
points which coordinates     are not expressible in the chosen
numeral system (recall Postulate~2). Thus, we are able to
calculate the number of points if we have already decided which
numerals will be used to express the coordinates of points.

Different numeral systems can be chosen to express coordinates of
the points in dependence  on the precision level we want to
obtain. For example, Pirah\~{a} are not able to express any point.
If the numbers $0 \le x < 1$ are expressed in the form
$\frac{p-1}{\mbox{\ding{172}}}, p \in \mathbb{N}$, then the
smallest positive number we can distinguish is
$\frac{1}{\mbox{\ding{172}}}$ and the interval $[0,1)$ contains
the following  points
 \beq
  0, \,\,\,\,
\frac{1}{\mbox{\ding{172}}}, \,\,\,\, \frac{2}{\mbox{\ding{172}}},
\,\,\,\,  \ldots \,\,\,\,
\frac{\mbox{\ding{172}}-2}{\mbox{\ding{172}}}, \,\,\,\,
\frac{\mbox{\ding{172}}-1}{\mbox{\ding{172}}}.
 \label{der2}
       \eeq
It is easy to see that they are \ding{172}. If  we want to count
the number of intervals of the form $[a-1,a), a \in \mathbb{N},$
on the ray  $x \ge 0$,  then, due to Postulate~3, the definition
of sequence, and Theorem~\ref{t2}, not more than \ding{172}
intervals of this type can be distinguished  on the ray  $x \ge
0$. They are
\[
[0,1), \,[1,2), \,[2,3),\, \ldots \,
 [\mbox{\ding{172}}-3,\mbox{\ding{172}}-2),\,[\mbox{\ding{172}}-2,\mbox{\ding{172}}-1),\,[\mbox{\ding{172}}-1,\mbox{\ding{172}}).
\]
Within each of them we are able to distinguish $\mbox{\ding{172}}$
points and, therefore, at the entire ray $\mbox{\ding{172}}^{2}$
 points can be observed. Analogously, the ray $x < 0$ is represented by the
intervals
\[
[-\mbox{\ding{172}},-\mbox{\ding{172}}+1), \,
[-\mbox{\ding{172}}+1,-\mbox{\ding{172}}+2),\,
 \ldots \,
 [-2,-1),\,[-1,0).
\]
Hence, this ray also contains $\mbox{\ding{172}}^{2}$ such points
and on the whole line   $2\mbox{\ding{172}}^{2}$ points of this
type can be represented and observed.

Note that the point $-\mbox{\ding{172}}$ is included in this
representation and the point $\mbox{\ding{172}}$ is excluded from
it. Let us slightly modify our numeral system in order to have
$\mbox{\ding{172}}$ representable. For this purpose, intervals of
the type $(a-1,a], a \in \mathbb{N},$ should be considered to
represent the ray $x > 0$ and the separate symbol, 0, should be
used to represent zero. Then, on the ray $x > 0$ we are able to
observe $\mbox{\ding{172}}^{2}$ points and, analogously, on the
ray $x < 0$ we also are able to observe $\mbox{\ding{172}}^{2}$
points. Finally, by adding the symbol used to represent zero we
obtain that on the entire line $2\mbox{\ding{172}}^{2}+1$ points
can be observed.

It is important to stress that the situation with counting points
is a direct consequence of Postulate~2 and is   typical for the
natural sciences where it is well known that instruments influence
results of observations. It is similar to the work with
microscopes or fractals (see \cite{fractals}): we decide the level
of the precision we need and obtain a result dependent on the
chosen level of accuracy. If we need a more precise or a more
rough answer, we change the lens  of our microscope.

 In our terms this means to
change one numeral system with another. For instance, instead of
the numerals considered above, let us choose a positional numeral
system with the radix\index{radix} $b$
 \beq
  (.a_{1} a_{2}
\ldots a_{q-1} a_{q})_b,  \hspace{1cm} q \in \mathbb{N},
  \label{der3}
       \eeq
to calculate the number of points within the interval $[0,1)$.

\begin{theorem}
\label{tder1} The number of elements   of the set
 of numerals of the type (\ref{der3}) is equal to
$b^{\mbox{\ding{172}}}$.
\end{theorem}

\textit{Proof.} Formula (\ref{der3}) defining the type of numerals
we deal with  contains a  sequence  of digits    $a_{1} a_{2}
\ldots a_{q-1} a_{q}$. Due to the definition of the sequence and
Theorem~\ref{t2}, this sequence can have as a maximum \ding{172}
elements, i.e., $q \le \mbox{\ding{172}}$. Thus, it can be at
maximum \ding{172} positions on the
  the right of the dot. Every position can
be filled in by one of the $b$ digits from the
alphabet\index{alphabet} $\{ 0, 1, \ldots , b-1 \}$. Thus, we have
$b^{\mbox{\ding{172}}}$ combinations. As a result,   the
positional numeral system using the numerals of the form
(\ref{der3}) can express $b^{\mbox{\ding{172}}}$ numbers. \hfill
 $\Box$
\begin{corollary}
\label{c1} The entire line contains
$2\mbox{\ding{172}}b^{\mbox{\ding{172}}}$ points of the type
(\ref{der3}).
\end{corollary}

\textit{Proof.}
 We have already
seen above that it is possible to distinguish $2\mbox{\ding{172}}$
unit intervals within the line. Thus,   the whole number of points
of the type (\ref{der3}) on the line is equal to
$2\mbox{\ding{172}}b^{\mbox{\ding{172}}}$. \hfill
 $\Box$

In this example of counting, we have changed the tool to calculate
the number of points within each unit interval from (\ref{der2})
to (\ref{der3}), but used the old way to calculate the number of
intervals, i.e., by natural numbers. If we are not interested in
subdividing the line at intervals and want to obtain the number of
the points on the line directly by using positional numerals of
the type
 \beq
 (a_{n-1}a_{n-2} \ldots a_1 a_0 .
 a_{1} a_{2}
\ldots a_{q-1} a_{q})_b
   \label{3.103}
       \eeq
 with $n, \, q \in \mathbb{N}$, then the following result holds.
\begin{corollary}
\label{c2} The number of elements of the set,$\mathbb{R}_{b}$, of
numerals of the type (\ref{3.103}) is $|\mathbb{R}_{b}| =
b^{2\mbox{\ding{172}}}$.
\end{corollary}

\textit{Proof.} In formula (\ref{3.103}) defining the type of
numerals we deal with there are two sequences of digits: the first
one, $a_{n-1}, a_{n-2}, \ldots a_1, a_0$, is used to express the
integer part of the number and the second, $a_{1}, a_{2}, \ldots
a_{q-1}, a_{q}$, for its fractional part. Analogously to the proof
of Theorem~\ref{tder1}, we can have as a maximum
$b^{\mbox{\ding{172}}}$ combinations to express the integer part
of the number and the same quantity to express its fractional
part. As a result,   the positional numeral system using the
numerals of the form (\ref{3.103})   can express
$b^{2\mbox{\ding{172}}}$ numbers. \hfill
 $\Box$

It is worth     noticing that   in our approach, all the numerals
from (\ref{3.103})   represent different numbers. It is possible
to execute, for example, the following subtraction
\[
3.\underbrace{0000\ldots000}_{\mbox{\ding{172}
positions}}-2.\underbrace{9999\ldots999}_{\mbox{\ding{172}
positions}} =0.\underbrace{0000\ldots001}_{\mbox{\ding{172}
positions}}.
\]
On the other hand, the traditional point of view on real numbers
tells us that there exist real numbers that can be represented in
positional systems by two different infinite sequences of digits,
for instance, in the decimal positional system the records
$3.000000\ldots$ and $2.99999\ldots$ represent the same number.
Note that there is no any contradiction between the traditional
and the new points of view. They just use different lens in their
mathematical microscopes to observe numbers. The instrument  used
in the traditional point of view for this purpose was just too
weak to distinguish two different numbers in the records
$3.000000\ldots$ and $2.99999\ldots$.

 We conclude this section by
the following observation. Traditionally, it was accepted that
\textit{any} positional numeral system is able to represent
\textit{all} real numbers (`the whole real line'). In this
section, we have shown that any numeral system is just an
instrument that can be used to observe \textit{certain} real
numbers. This instrument can be more or less powerful, e.g., the
positional system (\ref{3.103}) with the radix 10 is more powerful
than the positional system (\ref{3.103}) with the radix 2 but
neither of the two is able to represent irrational numbers (see
\cite{informatica}). Two numeral systems can allow us to observe
either the same sets of numbers, or sets of numbers having an
intersection, or two disjoint sets of numbers. Due to Postulate~2,
we are not able to answer the question `What is the whole real
line?' because this is the question asking `What is the object of
the observation?', we are able just to invent more and more
powerful numeral systems that will allow us to improve our
observations of numbers by using newly introduced numerals.

\section{Two concepts of continuity}
\label{s4}

The goal of   this section is to    develop a new point of view on
the notion of continuity that  would be, on the one hand, more
physical and, on the other hand,   could be used for functions
assuming infinite and infinitesimal values.

In Physics, the `continuity'\label{p:continuity}  of an object is
relative. For example, if we observe a table by eye, then we see
it as being continuous. If we use a microscope for our
observation, we see that the table is discrete. This means that
\textit{we decide} how to see the object, as a continuous or as a
discrete, by the choice of the instrument for observation. A weak
instrument -- our eyes -- is not able to distinguish its internal
small separate parts (e.g., molecules) and we see the table as a
continuous object. A sufficiently strong microscope allows us to
see the separate parts and the table becomes discrete but each
small part now is viewed as continuous.

In contrast, in   traditional mathematics, any mathematical object
is either continuous or discrete. For example, the same function
cannot be   both  continuous and discrete. Thus, this
contraposition of discrete and continuous in the traditional
mathematics does not reflect  properly the physical situation that
we observe in practice. The infinite and infinitesimal numbers
described in the previous sections give us a possibility to
develop a new theory of continuity that is closer to the physical
world and better reflects the new discoveries made by physicists.
Recall that the foundations of the mathematical analysis have been
established centuries ago and, therefore, do not take into account
the subsequent revolutionary results in Physics, e.g., appearance
of Quantum Physics (the goal of non-standard analysis was to
re-write these foundations by using non-archimedean ordered field
extensions of the reals and not to include Physics in Analysis).
In this section, we start by introducing a definition of the
one-dimensional continuous set of points based on Postulate~2  and
the above consideration and by establishing relations to such a
fundamental notion as a function using the infinite and
infinitesimal numbers.

We recall that traditionally a function   $f(x)$ is defined as a
binary relation among two sets $X$ and $Y$   (called the
\textit{domain} and the \textit{codomain} of the relation) with
the additional property that to each element $x \in X$ corresponds
exactly one element $f(x) \in Y$\label{p:function}. We now
consider a function $f(x)$ defined over a one-dimensional interval
$[a,b]$. It follows immediately from the previous sections  that
to define a function $f(x)$ over an interval $[a,b]$ it is not
sufficient to give a rule for evaluating $f(x)$ and the values $a$
and $b$ because we are not able to evaluate $f(x)$ at \textit{any}
point $x \in [a,b]$ (for example, traditional numeral systems do
not allow us to express \textit{any} irrational number $\zeta$
and, therefore, we are not able to evaluate $f(\zeta)$). However,
the traditional definition of a function includes in its domain
points at which $f(x)$ cannot be evaluated, thus introducing
ambiguity.

Note that a numeral system can include  certain numerals some of
which can be expressed as a  result of arithmetical operations
with other symbols    and some of them cannot. Such symbols as $e,
\pi, \sqrt{3},$ and other special symbols used to represent
certain irrational numbers cannot be expressed by any known
numeral system that uses only symbols representing  integer
numbers. These symbols are introduced in the mathematical language
by their properties as numerals 0 or 1. The introduction of
numerals $e, \pi, \sqrt{3},$ etc. in a numeral system,  of course,
enlarges its possibilities to represent numbers but, in any way,
these possibilities remain limited.

Thus, in order to be precise in the definition of a function, it
is necessary to indicate explicitly a numeral system,
$\mathcal{S}$, we intend to use to express points from the
interval $[a,b]$. A function $f(x)$ becomes defined when we know a
rule allowing us to obtain $f(x)$, given $x$ and its domain, i.e.,
the set $[a,b]_{\mathcal{S}}$ of points $x \in [a,b]$ expressible
in the chosen numeral system $\mathcal{S}$. We suppose hereinafter
that the   system $\mathcal{S}$  is used to write down $f(x)$ (of
course, the choice of   $\mathcal{S}$ determines a class of
formulae and/or procedures we are able to express using
$\mathcal{S}$) and it allows us to express any number
\[
y=f(x), \hspace{5mm} x \in [a,b]_{\mathcal{S}}.
\]
 The number of points of the domain
$[a,b]_{\mathcal{S}}$ can be finite or infinite but the set
$[a,b]_{\mathcal{S}}$ is always discrete. This means that for any
point $x \in [a,b]_{\mathcal{S}}$ it is possible to determine its
closest right and left neighbors, $x^{+}$ and $x^{-}$,
respectively, as follows
 \beq
x^{+}= \min \{z: z \in [a,b]_{\mathcal{S}},\hspace{3mm} z >x \},
\hspace{5mm}
 x^{-}=  \max \{z: z \in [a,b]_{\mathcal{S}}, \hspace{3mm} z <x \}.
\label{5.1}
 \eeq

Apparently, the obtained discrete construction leads us to the
necessity of abandoning  the nice idea of continuity, which is a
very useful notion used in different fields of mathematics. But
this is not the case. In contrast, the new approach allows us to
introduce a new definition of continuity   very well reflecting
the physical world.

Let us consider    $n+1$ points at a line
 \beq
a=x_{0} < x_{1} < x_{2} < \ldots < x_{n-1} < x_{n}=b
 \label{5.3}
 \eeq
and suppose that we have a numeral system  $\mathcal{S}$ allowing
us to calculate their coordinates   using a unit of measure $\mu$
(for example, meter, inch, etc.)  and to thus construct   the set
$X=[a,b]_{\mathcal{S}}$ expressing these points.

The set $X$ is called \textit{continuous in the unit of measure}
$\mu$ if for any $ x \in (a,b)_{\mathcal{S}}$ it follows that the
differences $x^{+}-x$ and $x-x^{-}$ from (\ref{5.1}) expressed in
units $\mu$ are equal to infinitesimal numbers. In our numeral
system with radix grossone this means that all the differences
$x^{+}-x$ and $x-x^{-}$ contain only negative grosspowers. Note
  that it becomes possible to differentiate types of continuity
by taking into account values of grosspowers of infinitesimal
numbers (continuity of order $\mbox{\ding{172}}^{-1},$ continuity
of order $\mbox{\ding{172}}^{-2}$, etc.).

 \begin{figure}[t]
  \begin{center}
    \epsfig{ figure = 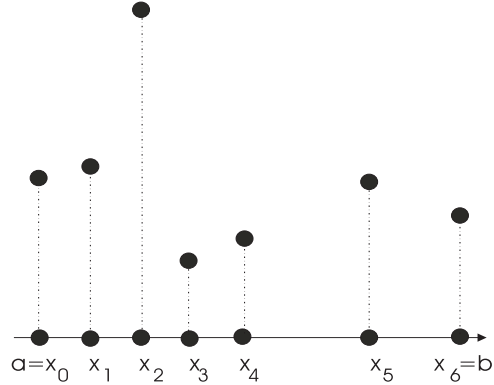, width = 2.5in, height = 2in,  silent = yes }
    \caption{It is not possible to say whether this function is continuous or discrete until
    we have not introduced a unit of measure and a numeral system to express distances
    between the points}
 \label{Big_paper4}
  \end{center}
\end{figure}

This  definition emphasizes the physical principle that there does
not exist an absolute continuity: it is relative  with respect to
the chosen instrument of observation which  in our case is
represented by the unit of measure~$\mu$. Thus, the same set can
be viewed as a continuous or not depending on the chosen unit of
measure.

\begin{example}
\label{e13} The set of five equidistant points
 \beq
X_1 = \{ a,   x_{1}, x_{2}, x_{3}, x_{4}  \}
 \label{5.3.0}
 \eeq
from Fig.~\ref{Big_paper4} can have the distance $d$ between the
points equal to $\mbox{\ding{172}}^{-1}$ in a unit of measure
$\mu$ and to be, therefore, continuous in $\mu$. Usage of a new
unit of measure $\nu = \mbox{\ding{172}}^{-3}\mu$ implies that
$d=\mbox{\ding{172}}^{2}$ in $\nu$ and the set $X_1$ is not
continuous in $\nu$. \hfill $\Box$
\end{example}

Note that the introduced definition does not require that all the
points from $X$ are equidistant. For instance, if in
Fig.~\ref{Big_paper4} for a unit measure $\mu$ the largest over
the set $[a,b]_{\mathcal{S}}$ distance $x_{5}-x_{4}$ is
infinitesimal, then the whole set  is continuous in $\mu$.

The set $X$ is called \textit{discrete in the unit of measure}
$\mu$ if for all points $ x \in (a,b)_{\mathcal{S}}$ it follows
that the differences $x^{+}-x$ and $x-x^{-}$ from (\ref{5.1})
expressed in units $\mu$ are not infinitesimal numbers. In our
numeral system with radix grossone this means that in all the
differences $x^{+}-x$ and $x-x^{-}$   negative grosspowers cannot
be the  largest ones. For instance, the set $X_1$ from
(\ref{5.3.0}) is discrete in the unit of measure $\nu$ from
Example~\ref{e13}. Of course, it is also possible to consider
intermediate cases where sets have continuous and discrete parts.

The introduced notions allow  us to give the following very simple
definition of a function continuous at a point. A function $f(x)$
defined over a   set $[a,b]_{\mathcal{S}}$  continuous in a unit
of measure  $\mu$ is called \textit{continuous in the unit of
measure  $\mu$ at a point} $x \in (a,b)_{\mathcal{S}}$ if both
differences $f(x)-f(x^{+})$ and $f(x)-f(x^{-})$ are infinitesimal
numbers in   $\mu$, where $x^{+}$ and $x^{-}$ are from
(\ref{5.1}). For the continuity at points $a$, $b$ it is
sufficient that one of these differences be infinitesimal. The
notions of continuity from the left and from the right in a unit
of measure $\mu$ at a point are introduced naturally. Similarly,
the notions of a function discrete, discrete from the right, and
discrete from the left can be defined\footnote{Note that in these
definitions we have accepted that the same unit of measure, $\mu$,
has been used to measure distances along both axes, $x$ and
$f(x)$. A natural generalization can be done in case of need by
introducing two different units of measure, let say, $\mu_1$ and
$\mu_2$, for measuring distances along two axes.}.

The function $f(x)$ is \textit{continuous in the unit of measure
$\mu$ over the set} $[a,b]_{\mathcal{S}}$ if it is continuous in
$\mu$ at all points of $[a,b]_{\mathcal{S}}$. Again, it becomes
possible to differentiate types of continuity by taking into
account values of grosspowers of infinitesimal numbers (continuity
of order $\mbox{\ding{172}}^{-1},$ continuity of order
$\mbox{\ding{172}}^{-2}$, etc.) and to consider functions in such
units of measure that they become continuous or discrete over
certain subintervals of $[a,b]$.
 Hereinafter, we shall often fix the unit of measure $\mu$ and
write just `continuous function' instead of `continuous function
in the unit of measure  $\mu$'. Let us give some examples
illustrating the introduced definitions.

\begin{example}
\label{e14}
 We start by showing that the function $f(x)=x^2$ is
continuous over the set $X_2$\label{p:X2} defined as the interval
$[0,1]$ where numerals $\frac{i}{\mbox{\ding{172}}}, 0 \le i \le
\mbox{\ding{172}},$ are used to express its points in units $\mu$.
First of all, note that the set $X_2$ is continuous in   $\mu$
because its points are equidistant with the distance
$d=\mbox{\ding{172}}^{-1}$. Since this function is strictly
increasing,  to show its continuity it is sufficient to check the
difference $f(x)-f(x^{-})$ at the point $x=1$. In this case,
$x^{-}=1-\mbox{\ding{172}}^{-1}$ and we have
\[
 f(1)-f(1-\mbox{\ding{172}}^{-1})=
 1-(1-\mbox{\ding{172}}^{-1})^2 =
 2\mbox{\ding{172}}^{-1}(-1)\mbox{\ding{172}}^{-2}.
\]
This number is infinitesimal, thus $f(x)=x^2$ is continuous over
the set $X_2$. \hfill $\Box$
\end{example}

\begin{example}
\label{e15} Consider the same function $f(x)=x^2$ over the set
$X_3$ defined as the interval
$[\mbox{\ding{172}}-1,\mbox{\ding{172}}]$ where numerals
$\mbox{\ding{172}}-1+\frac{i}{\mbox{\ding{172}}}, 0 \le i \le
\mbox{\ding{172}},$ are used to express its points  in units
$\mu$. Analogously, the set $X_3$ is continuous and it is
sufficient to check the difference $f(x)-f(x^{-})$ at the point
$x=\mbox{\ding{172}}$ to show continuity of $f(x)$ over this set.
In this case,
\[
x^{-}=\mbox{\ding{172}}-1+\frac{\mbox{\ding{172}}-1}{\mbox{\ding{172}}}=
\mbox{\ding{172}}-\mbox{\ding{172}}^{-1},
\]
\[
f(x)-f(x^{-}) =
f(\mbox{\ding{172}})-f(\mbox{\ding{172}}-\mbox{\ding{172}}^{-1})=
\mbox{\ding{172}}^2-(\mbox{\ding{172}}-\mbox{\ding{172}}^{-1})^2
=2\mbox{\ding{172}}^{0}(-1)\mbox{\ding{172}}^{-2}.
\]
This number is not infinitesimal because it contains the finite
part $2\mbox{\ding{172}}^{0}$ and, as a consequence, $f(x)=x^2$ is
not continuous over the set $X_3$. \hfill $\Box$
\end{example}

\begin{example}
\label{e16} Consider   $f(x)=x^2$   defined over the set $X_4$
being the interval $[\mbox{\ding{172}}-1,\mbox{\ding{172}}]$ where
numerals $\mbox{\ding{172}}-1+\frac{i}{\mbox{\ding{172}}^2}, 0 \le
i \le \mbox{\ding{172}}^2,$ are used to express its points   in
units $\mu$. The set $X_4$ is continuous and we check the
difference $f(x)-f(x^{-})$ at the point $x=\mbox{\ding{172}}$. We
have
\[
x^{-}=\mbox{\ding{172}}-1+\frac{\mbox{\ding{172}}^2-1}{\mbox{\ding{172}}^2}=
\mbox{\ding{172}}-\mbox{\ding{172}}^{-2},
\]
\[
f(x)-f(x^{-}) =
f(\mbox{\ding{172}})-f(\mbox{\ding{172}}-\mbox{\ding{172}}^{-2})=
\mbox{\ding{172}}^2-(\mbox{\ding{172}}-\mbox{\ding{172}}^{-2})^2
=2\mbox{\ding{172}}^{-1}(-1)\mbox{\ding{172}}^{-4}.
\]
Since the obtained result is   infinitesimal, $f(x)=x^2$ is
continuous over   $X_4$. \hfill $\Box$
\end{example}

Let us consider now a function $f(x)$ defined by formulae over a
set $X=[a,b]_{\mathcal{S}}$ so that different expressions can be
used over different subintervals of $[a,b]$. The term `formula'
hereinafter indicates a single expression being a sequence of
numerals and arithmetical operations used to evaluate
$f(x)$\label{p:formula}.

\begin{example}
\label{e17}    The function $g(x)=2x^2-1,     x \in
[a,b]_{\mathcal{S}},$   is defined by one formula and   function
 \beq
f(x) = \left\{ \begin{array}{ll} \max \{ -14x, 25x^{-1}  \}, &   x
\in
[c,0)_{\mathcal{S}} \cup (0,d]_{\mathcal{S}},\\
2x, & x=0, \end{array} \right.  \hspace{5mm} c <0,  \hspace{5mm} d
> 0,
 \label{5.3.1-1}
 \eeq
is defined by three formulae, $f_1(x), f_2(x),$ and $f_3(x)$ where
 \beq
\begin{array}{ll}
f_1(x) =     -14x,  &   x \in [c,0)_{\mathcal{S}},\\
f_2(x) =     2x,  &   x =0,\\
f_3(x) =      25x^{-1},  &   x \in (0,d]_{\mathcal{S}}.
\hspace{5mm}
 \Box
\end{array}
 \label{5.3.1}
 \eeq
\end{example}
\noindent Consider now a function $f(x)$ defined in a neighborhood
of a point $x$ as follows
 \beq
f(\xi) = \left\{ \begin{array}{ll}
f_1(\xi),  &   x-l \le  \xi < x,\\
f_2(\xi),  &   \xi = x,\\
f_3(\xi),  &   x<\xi \le x+r,  \end{array}
 \right.
 \label{5.3.2}
 \eeq
where the number $l$ is any number such that the same formula
$f_1(\xi)$ is used to define $f(\xi)$ at all points  $\xi$ such
that $x-l \le  \xi < x$. Analogously, the number $r$ is any number
such that the same formula $f_3(\xi)$ is used to define $f(\xi)$
at all points $\xi$ such that $x<\xi \le x+r$.  Of course, as a
particular case it is possible that the same formula is used to
define $f(\xi)$ over the interval $[x-l,x+r]$, i.e.,
  \beq
f(\xi)=f_1(\xi)=f_2(\xi)=f_3(\xi),  \hspace{5mm} \xi \in
[x-l,x+r].
 \label{5.3.3}
 \eeq
It is also possible that (\ref{5.3.3}) does not hold but formulae
$f_1(\xi)$ and $f_3(\xi)$ are defined at the point $x$ and are
such that at this point they return the same value, i.e.,
  \beq
f_1(x)=f_2(x)=f_3(x).
 \label{5.3.3.1}
 \eeq
If condition (\ref{5.3.3.1}) holds,   we say that function $f(x)$
has \textit{continuous formulae} at the point $x$. Of course,  in
the general case, formulae $f_1(\xi),f_2(\xi),$ and $f_3(\xi)$ can
be or cannot be defined out of the respective intervals from
(\ref{5.3.2}).  In cases where condition (\ref{5.3.3.1}) is not
satisfied we say that function $f(x)$ has \textit{discontinuous
formulae}   at the point $x$. Definitions of functions having
formulae which are continuous or discontinuous from the left and
from the right are introduced naturally. Let us give an example
showing that the introduced definition can also be easily used for
functions assuming infinitesimal and infinite values.
\begin{example}
\label{e23} Let us study the following function that at infinity,
in the neighborhood of the point $x=\mbox{\ding{172}}$, assumes
infinitesimal values
 \beq
f(x) = \left\{ \begin{array}{ll}
\mbox{\ding{172}}^{-1} + 5(x-\mbox{\ding{172}}),  &   x \neq \mbox{\ding{172}},\\
\mbox{\ding{172}}^{-2},  &   x=\mbox{\ding{172}}.  \end{array}
 \right.
  \label{der6}
 \eeq
 By using designations (\ref{5.3.2})   we have
 \[ f(\xi) = \left\{ \begin{array}{ll}
f_1(\xi)=\mbox{\ding{172}}^{-1} + 5(x-\mbox{\ding{172}}),  &    \xi < \mbox{\ding{172}},\\
f_2(\xi)=\mbox{\ding{172}}^{-2},  &   \xi = \mbox{\ding{172}},\\
f_3(\xi)=\mbox{\ding{172}}^{-1} + 5(x-\mbox{\ding{172}}),  & \xi >
\mbox{\ding{172}},  \end{array}
 \right.
 \]
Since
 \[
f_1(\mbox{\ding{172}})=f_3(\mbox{\ding{172}})=
\mbox{\ding{172}}^{-1} \neq
f_2(\mbox{\ding{172}})=\mbox{\ding{172}}^{-2},
 \]
we conclude that     the function (\ref{der6}) has discontinuous
formulae at the point $x=\mbox{\ding{172}}$.  It is remarkable
that we were able to establish this  easily in spite of the fact
that all three values, $f_1(\mbox{\ding{172}}),
f_2(\mbox{\ding{172}}),$ and $ f_3(\mbox{\ding{172}})$  were
infinitesimal and were evaluated at infinity.  Analogously, the
function (\ref{5.3.1-1}) has continuous formulae at the point
$x=0$ from the left and discontinuous from the right.
 \hfill$\Box$
\end{example}

\begin{example}
\label{e18} Let us study the following function
 \beq
f(x) = \left\{ \begin{array}{ll}
\mbox{\ding{172}}^3+\frac{x^2-1}{x-1},  &   x \neq 1,\\
a,  &   x=1,  \end{array}
 \right.
 \label{5.4}
 \eeq
at the point $x=1$. By using designations (\ref{5.3.2}) and the
fact that for $x \neq 1$ it follows $\frac{x^2-1}{x-1}=x+1$, from
where we have
 \[ f(\xi) = \left\{ \begin{array}{ll}
f_1(\xi)=\mbox{\ding{172}}^3+\xi+1,  &    \xi < 1,\\
f_2(\xi)=a,  &   \xi = 1,\\
f_3(\xi)=\mbox{\ding{172}}^3+\xi+1,  &    \xi > 1,  \end{array}
 \right.
 \]
Since
 \[ f_1(1)=f_3(1)=\mbox{\ding{172}}^3+2, \hspace{1cm} f_2(1)=a,
 \]
we obtain that if $a=\mbox{\ding{172}}^3+2$,  then the function
(\ref{5.4}) has   continuous formulae at the point $x=1$,
otherwise it has discontinuous formulae at this point.
 Note, that even if
$a=\mbox{\ding{172}}^3+2+\varepsilon$, where $\varepsilon$ is an
infinitesimal number (remind that all infinitesimals are not equal
to zero), we   establish easily that the function has
discontinuous formulae in spite of the fact that both numbers,
$\mbox{\ding{172}}^3+2$ and $a$, are infinite.
 \hfill$\Box$
\end{example}

Similarly to the existence of numerals that cannot be expressed
through other numerals, in mathematics there exist   functions
that cannot be expressed as  a sequence of numerals and
arithmetical operations because (again similarly to numerals) they
are introduced through their properties. Let us consider, for
example, the   function $f(x)=\sin(x)$. It can be immediately seen
that it does not satisfy the traditional definition of a function
(see the corresponding discussion in page~\pageref{p:function}).
We know well its codomain but the story becomes more difficult
with its domain and impossible with the relation that it is
necessary to establish to obtain $f(x)$ when a value for $x$ is
given. In fact, we know precisely the value of $\sin(x)$ only at
certain points $x$; for other points the value of $\sin(x)$ are
measured (as a result, errors are introduced) or approximated
(again errors are introduced), and, finally, as it has been
already mentioned, not all points $x$ can be expressed by known
numeral systems.

However, traditional approaches allow us to confirm its continuity
(which is clear due to the physical way it has been introduced)
also from positions of   general definitions of continuity both in
standard and non-standard Analysis. Let us show that, in spite of
the fact that we do not know a complete formula for calculating
$\sin(x)$, the new approach also allows us to   show that
$\sin(x)$ is a function with continuous formulae. This can be done
by using   partial information about the structure of the formula
of $\sin(x)$ by appealing to the same geometrical ideas that are
used in the traditional proof of continuity of $\sin(x)$
(cf.~\cite{Courant}).

\begin{theorem}
\label{tder2} The function $f(x)=\sin(x)$ has   continuous
formulae.
 \end{theorem}

\textit{Proof.} Let us first show that $\sin(x)$ has   continuous
formulae at $x=0$. By using designations (\ref{5.3.3}) and
well-known geometrical considerations  we can write for $f_1(\xi),
f_3(\xi)$ and $x=0$ the following relations
 \beq
 \xi < f_1(\xi) < 0, \hspace{5mm} -\frac{\pi}{2} < \xi < 0,
 \label{der4}
\eeq
 \beq
 0 < f_3(\xi) < \xi, \hspace{5mm} 0 < \xi < \frac{\pi}{2}.
 \label{der5}
\eeq
 Thus, even though we are not able to calculate $f_1(\xi)$ and
$f_3(\xi)$ at points $\xi \neq 0$, we can use the estimates
(\ref{der4}) and (\ref{der5}) from where we obtain
$f_1(0)=f_3(0)=0$ because the estimate $\xi$ being at the left
part of (\ref{der4}) and at the right part of (\ref{der5}) can be
evaluated at the point  $\xi = 0$. It then follows from the
obvious fact $f_2(\xi)= 0$ that $\sin(x)$ has   continuous
formulae at $x=0$. By a complete analogy, the fact that the
function $f(x)=\cos(x)$ has continuous formulae at $x=0$ can be
proved.

If we represent now the point $\xi$ in the form $\xi=x+\zeta$ then
we can write
\[
\sin(x+\zeta)=\sin(x)\cos(\zeta)+\cos(x)\sin(\zeta).
\]
Inasmuch as both $\sin(x)$ and $\cos(x)$ have continuous formulae
at the point $x=0$, it follows that at the point $\xi=x$ we have
\[
f_1(x)=f_3(x)=\sin(x+0)=\sin(x)\cos(0)+\cos(x)\sin(0)=\sin(x).
\]
This fact concludes the proof because obviously $ f_2(x)=\sin(x)$,
too. \hfill
 $\Box$

 \begin{corollary}
\label{c3} The function $f(x)=\cos(x)$ has   continuous formulae.
\end{corollary}
\textit{Proof.} The fact is a straightforward consequence of the
theorem.  \hfill
 $\Box$

To conclude this section we emphasize that functions having
continuous formulae at a point can be continuous or discrete at
this point depending on the chosen unit of measure. Analogously,
functions having discontinuous formulae at a point can be
continuous or discrete at this point again depending on the chosen
unit of measure.  The notion of continuity of a function depends
on the chosen unit of measure and numeral system $\mathcal{S}$ and
it can be used for functions defined by formulae, computer
procedures, tables, etc. In contrast, the notion of a function
having continuous formulae works only for functions defined by
formulae and does not depend on units of measure or numeral
systems  chosen to express its domain.  It is related only to
properties of formulae and does not depend on the domain at all.
Note that we have established the facts of continuity or
discontinuity of formulae in Theorem~\ref{tder2} and
Examples~\ref{e23}, \ref{e18}  without indicating domains of the
considered functions.

\section{Differential calculus}
\label{s5}

In this section, the notions of the first derivative and
subdifferential are introduced. Both concepts are defined without
usage of limits and for functions assuming finite, infinite, and
infinitesimal values. We first give all the definitions and then
illustrate them by a series of examples. Special attention (as in
the entire paper) is paid to the computational issues and their
relations to the introduced definitions.

We shall call the following two expressions, $f^{-}(x,l)$ and
$f^{+}(x,r)$, that are obtained by isolating the multipliers $l$
and $r$ in the left and the right differences $f_1(x)-f_1(x-l)$
and $f_3(x+r)-f_3(x)$, respectively, \textit{the left and the
right relative differences for} $f(x)$ \textit{at a point} $x$:
 \beq
f_1(x)-f_1(x-l)=  l \cdot f^{-}(x,l), \hspace{1cm}   l \ge 0,
 \label{5.3.4}
 \eeq
 \beq
f_3(x+r)-f_3(x)  =  r \cdot f^{+}(x,r), \hspace{1cm}   r \ge 0.
 \label{5.3.5}
 \eeq
In the introduced definitions (examples will be given soon) it is
not required that the numbers $l$ and $r$ from (\ref{5.3.2}) tend
to zero -- they can be finite or even infinite. Their boundaries
are determined by formulae $f_1(\xi),f_2(\xi),$ and $f_3(\xi)$
defining $f(\xi)$. Note that the value of $f(\xi)$ at the point
$\xi=x$ is defined by formula $f_2(\xi)$ which is not used in
(\ref{5.3.4}) and (\ref{5.3.5}).

If  for $\xi=x$  formula   $f^{-}(\xi,l)$ is defined at $l=0$
and/or $f^{+}(\xi,r)$ is defined at $r=0$, then functions
 \beq
f^{-}(x)= f^{-}(x,0), \hspace{1cm} f^{+}(x) = f^{+}(x,0),
 \label{5.4.3.0}
 \eeq
are called  \textit{left} and \textit{right derivatives at the
point} $x$, respectively.   We can also introduce functions
 \beq
\tilde{f}^{-}(x,l) = f^{-}(x,l) - f^{-}(x),
 \label{5.4.3}
 \eeq
  \beq
\tilde{f}^{+}(x,r) = f^{+}(x,r) - f^{+}(x).
 \label{5.4.4}
 \eeq
Note that obtained formulae  give us possibility to evaluate
$f^{-}(x,l), \tilde{f}^{-}(x,l)$ and $f^{+}(x,r),
\tilde{f}^{+}(x,r)$  at any points
 \beq
x \in [a,b], \hspace{1cm}  x-l \in [a,b], \hspace{1cm}  x+r \in
[a,b],
 \label{5.4.1}
 \eeq
where formulae $f_1(\xi),f_2(\xi),$ and $f_3(\xi)$ are used but we
are able to calculate them only at  those points that can be
expressed in our numeral system $\mathcal{S}$, i.e., at the points
$x, x-l$, and $x+r$ such that
 \beq
x \in [a,b]_{\mathcal{S}}, \hspace{1cm} x-l \in
[a,b]_{\mathcal{S}}, \hspace{1cm} x+r \in [a,b]_{\mathcal{S}}.
 \label{5.4.2}
 \eeq

Geometrically,  $f^{-}(x,l)$ is the slope of the straight line
passing  through the points $(x-l,f_1(x-l))$ and $(x,f_1(x))$.
Analogously, $f^{+}(x,r)$ is the slope of the straight line
passing  through the points $(x,f_3(x))$ and $(x+r,f_3(x+r))$. The
left derivative, $f^{-}(x)$, can be viewed as the slope of the
line constructed at the point $(x,f_1(x))$ and
$\tilde{f}^{-}(x,l)$ is the difference between the slopes of this
line and that passing through the points $(x-l,f_1(x-l))$ and
$(x,f_1(x))$ with the slope $f^{-}(x,l)$. Analogously, $f^{+}(x)$
can be viewed as the slope of the line constructed at the point
$(x,f_3(x))$ and $\tilde{f}^{+}(x,r)$  is the difference between
the slopes of this line and that passing through the points
$(x,f_3(x))$ and $(x+r,f_3(x+r))$  with the slope $f^{+}(x,r)$.

Thus, the left derivative describes the behavior of $f(\xi)$ on
the left of the point $x$ and the right derivative on its right.
Both of them are independent of the value   $f_2(x)$ at the point
$x$.

Below we  introduce the notions  of derivative and derivatives
interval and give their geometrical interpretation in dependence
on the mutual positions of the two lines passing through the
points $(x,f_1(x))$ and $(x,f_3(x))$ with the slopes $f^{-}(x)$
and $f^{+}(x)$, respectively. We suppose that the right and the
left derivatives exist and function $f(x)$ has continuous formulae
at the point $x$. The latter supposition means that the points
$(x,f_1(x)), (x,f_2(x)),$ and $(x,f_3(x))$ coincide.

If at a point $x$ formulae $f^{-}(x)$ and $f^{+}(x)$ can be
written down in the same form then the following function
 \beq
f'(x)=f^{-}(x) = f^{+}(x)
 \label{5.4.5}
 \eeq
is called \textit{derivative} of $f(x)$ at the point $x$ and
function $f(x)$ is called \textit{antiderivative} of $f'(x)$ at
this point (this terminology will be used for $f(x)$ in the two
following cases too).

If at a point $\tau \in [a,b]_{\mathcal{S}}$ we obtain that
$f^{-}(\tau) = f^{+}(\tau)$ but formulae $f^{-}(x)$ and $f^{+}(x)$
cannot be written down in the same form, then the value
$f^{-}(\tau)$ is called \textit{derivative} of $f(x)$ at the point
$\tau$ and formulae to express derivative $f'(\tau)$ should be
chosen in concordance with the choice of the formula used to
express antiderivative $f(x)$. Geometrically, this case and the
previous one have the same meaning and $f'(\tau)$ is   the slope
of the straight line passing  through the points $(\tau,f(\tau))$
and tangent to the graph of  the function $f(x)$ at the point
$\tau$.

 \begin{figure}[t]
  \begin{center}
    \epsfig{ figure = 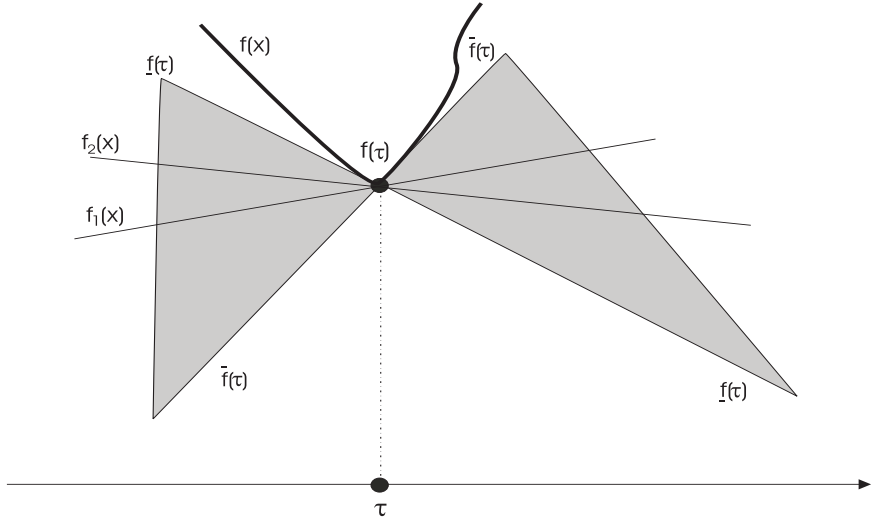, width =  4.5in, height = 2.7in,  silent = yes }
    \caption{Two tangent lines, $f_{1}(x)$ and $f_{2}(x)$, passing through the point $(\tau,f(\tau))$ and
    having slopes equal to   two derivatives
    $f'_{1}$ and $f'_{2}$ from the derivatives
    interval $[\underline{f}(\tau),\overline{f}(\tau)]$
     (shown in grey color) of a function $f(x)$ at a point $\tau$}
 \label{Big_paper5}
  \end{center}
\end{figure}

If at a point $\tau \in [a,b]_{\mathcal{S}}$ we obtain that
$f^{-}(\tau) \neq f^{+}(\tau)$ then the interval
$[\underline{f}(\tau),\overline{f}(\tau)]$ is called
\textit{derivatives interval of} $f(x)$ \textit{at the point}
$\tau$ where
 \beq
\underline{f}(\tau)= \min \{f^{-}(\tau), f^{+}(\tau) \},
 \label{5.4.6}
 \eeq
 \beq
\overline{f}(\tau) = \max \{f^{-}(\tau), f^{+}(\tau) \}.
 \label{5.4.7}
 \eeq
Geometrically, the derivatives interval contains slopes of all
lines tangent to the graph of the function at the point $\tau$. In
this case, the user chooses from the derivatives interval that
derivative which fits better his/her requirements (accuracy, type
of used algorithms, etc.) and works with it (a similar situation
takes place when one deals with subdifferentials (see
\cite{Clarke})). This case is illustrated in Fig.~\ref{Big_paper5}
were two tangent lines using derivatives, $f'_{1}$ and $f'_{2}$,
from the interval $[\underline{f}(\tau),\overline{f}(\tau)]$ are
shown.

Let us make a few comments. First of all, it is important to
notice that in the traditional approach the derivative (if it
exists) is a finite number and it can be defined only for
continuous functions. In our terminology, it can be finite,
infinite, or infinitesimal and its existence does not depend on
continuity of the function but on continuity of formulae. The
derivatives interval can exist  for discrete and continuous
functions defined over continuous sets, for functions defined over
discrete sets, and for functions having infinite or infinitesimal
values defined over sets having infinitesimal or infinite
boundaries.

The derivative (or derivatives interval) can be possibly defined
if the function under consideration has been introduced by
formulae and it cannot be defined if the function has been
introduced by a computer procedure or by a table. Thus, we
emphasize that just the fact of   presence of an analytical
expression of $f(x)$ allows us to find its derivatives interval.

When a function $f(x)$ is defined only by a computer procedure and
its analytical expressions are unknown, we cannot define
derivatives. We can only try to obtain a numerical approximation
without any estimate of its accuracy. For a  point $x \in
[a,b]_{\mathcal{S}}$ we define the interval
$[\underline{f}(x)_{\mathcal{S}},\overline{f}(x)_{\mathcal{S}}]$
called \textit{numerical derivatives interval of} $f(x)$
\textit{at the point} $x$ expressed in the system $\mathcal{S}$
where
 \beq
\underline{f}(x)_{\mathcal{S}}= \min \{f^{-}(x,x^{-}),
f^{+}(x,x^{+}) \},
 \label{5.5}
 \eeq
 \beq
\overline{f}(x)_{\mathcal{S}} = \max \{f^{-}(x,x^{-}),
f^{+}(x,x^{+}) \},
 \label{5.6}
 \eeq
and the points $x^{-}, x^{+}$ are from (\ref{5.1}). Each element
$f'(x) \in
[\underline{f}(x)_{\mathcal{S}},\overline{f}(x)_{\mathcal{S}}]$ is
called \textit{numerical derivative of} $f(x)$ \textit{at the
point} $x$ expressed in $\mathcal{S}$.

\begin{figure}[t]
  \begin{center}
    \epsfig{ figure = 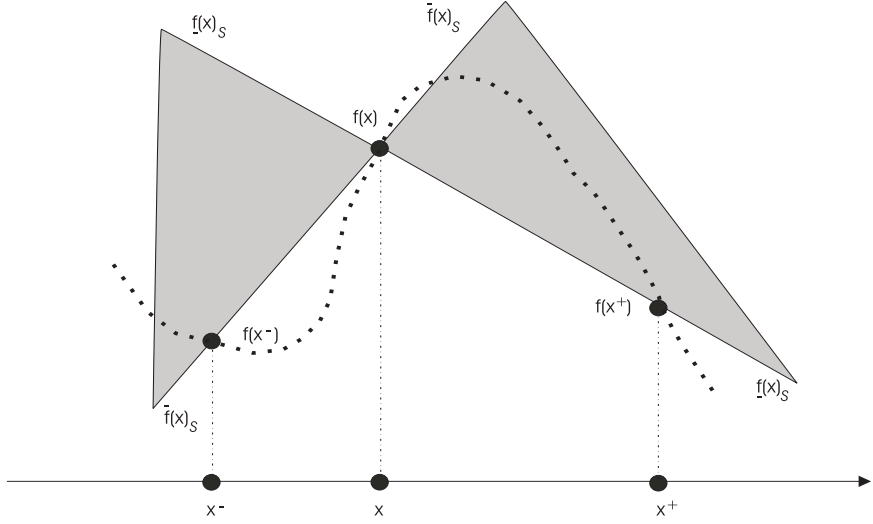, width = 4.5in, height = 2.7in,  silent = yes }
    \caption{The numeral system $\mathcal{S}$ able to express only points $x, x^{-},$ and
$x^{+}$ is too crude to present adequately the behavior of $f(x)$
} \label{Big_paper8}
  \end{center}
\end{figure}

Numerical derivatives interval can be useful also when formulae of
$f(x)$ are known. For example, if we have obtained that at a point
$x \in [a,b]_{\mathcal{S}}$ there exists the unique value $f'(x)$
and $f'(x) \notin
[\underline{f}(x)_{\mathcal{S}},\overline{f}(x)_{\mathcal{S}}]$
then this means that the chosen numeral system $\mathcal{S}$ is
too rude to be used for presentation of  $f(x)$. This situation is
illustrated in Fig.~\ref{Big_paper8} where  points $x, x^{-},$ and
$ x^{+}$ expressed in $\mathcal{S}$ are shown by big dots. Small
dots show the behavior of $f(x)$  expressed in a more powerful
numeral system. It can be clearly seen from Fig.~\ref{Big_paper8}
that $f'(x) \notin
[\underline{f}(x)_{\mathcal{S}},\overline{f}(x)_{\mathcal{S}}]$.

Note also that relaxed definitions of derivative and derivatives
interval can be also given by asking,  instead of (\ref{5.3.3.1}),
satisfaction of $f_1(\tau)=f_3(\tau)$  only (therefore, $f_1(\tau)
\neq f_2(\tau)$) can occur). For this case,   the same derivative
interval is obtained but it contains slopes of lines tangent to
the graph of the following function
\[
\hat{f}(x)=\left\{ \begin{array}{ll} f(x), & x \neq
 \tau,\\
 f_1(\tau),  & x = \tau.
 \end{array} \right.
\]

\begin{example}
\label{e19}

Consider again the function $f(x)=x^2$ at a point $x$  and
calculate its left and right differences. Condition (\ref{5.3.3})
holds for this function everywhere and, therefore, it has
continuous formulae at any point $x$. Then
\[
f_1(x)-f_1(x-l) = x^2-(x-l)^2  =    2xl-l^2   = l(2x-l)
 \]
and due to (\ref{5.3.4}) it follows
\[
 f^{-}(x,l) = 2x-l.
 \]
 Analogously, from (\ref{5.3.5}) we obtain
  \[
f_3(x+r)-f_3(x)  =  (x+r)^2 - x^2 =  2xr+r^{2}  = r(2x+r),
\]
  \[
f^{+}(x,r) =  2x+r.
\]
Since both functions $f^{-}(x)=f^{-}(x,0)$ and
$f^{+}(x)=f^{+}(x,0)$ exist and
 \[
 f^{-}(x) = f^{+}(x)=2x,
\]
the function $f(x)=x^2$ has the unique  derivative $f'(x)= 2x$ at
any point (finite, infinite, or infinitesimal) expressible in the
chosen numeral system $\mathcal{S}$. For instance, if
$x=\mbox{\ding{172}}$ then we obtain infinite values
$f(\mbox{\ding{172}})= \mbox{\ding{172}}^{2}$ and
$f'(\mbox{\ding{172}})= 2\mbox{\ding{172}}$. For infinitesimal
$x=\mbox{\ding{172}}^{-1}$ we have infinitesimal values
$f(\mbox{\ding{172}}^{-1})= \mbox{\ding{172}}^{-2}$ and
$f'(\mbox{\ding{172}}^{-1})= 2\mbox{\ding{172}}^{-1}$. \hfill
$\Box$
\end{example}
 It is important to emphasize that the introduced definitions
also work for functions including   infinite and infinitesimal
values in their formulae. For instance, it follows immediately
from the above consideration in Example~\ref{e19} that functions
$g_1(x)=\mbox{\ding{172}}^{-1}x^2$ and
$g_2(x)=\mbox{\ding{172}}x^2$ have derivatives
$g_1'(x)=2\mbox{\ding{172}}^{-1}x$ and
$g_2'(x)=2\mbox{\ding{172}}x$, respectively.

Inasmuch as the case where at a point $\tau \in
[a,b]_{\mathcal{S}}$ we obtain that $f^{-}(\tau) = f^{+}(\tau)$
but formulae $f^{-}(x)$ and $f^{+}(x)$ cannot be written down in
the same form is a simple consequence of the previous case, we
give an example of the situation where at a point $\tau \in
[a,b]_{\mathcal{S}}$ we obtain that $f^{-}(\tau) \neq
f^{+}(\tau)$.

\begin{example}
\label{e20}
 Let us consider the following function
 \beq
f(x)  = \left\{ \begin{array}{ll} -x, &
x \le 0,\\
 x+C,  & x >0,
 \end{array} \right.
 \label{5.6.1}
 \eeq
where $C$ is a constant. Then to define the function completely it
is necessary to choose an interval $[a,b]$ and   a numeral system
$\mathcal{S}$.

Suppose now that the chosen numeral system, $\mathcal{S}$, is such
that the point $x=0$ does not belong to $[a,b]_{\mathcal{S}}$.
Then  we immediately obtain from (\ref{5.3.4}) -- (\ref{5.4.3.0}),
(\ref{5.4.5}) that at any point $x \in
 [a,b]_{\mathcal{S}}$ the function has a unique derivative and
 \beq
f'(x) =  \left\{ \begin{array}{rl} -1, &
x < 0,\\
 1,  & x >0.
 \end{array} \right.
\label{5.6.2}
 \eeq

 \begin{figure}[t]
  \begin{center}
    \epsfig{ figure = 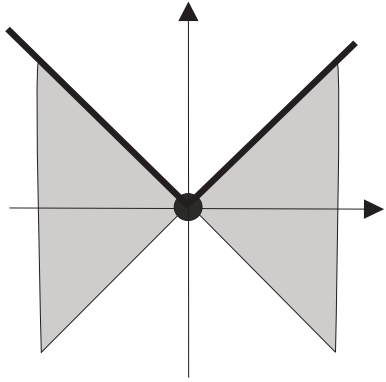, width = 2in, height = 2in,  silent = yes }
    \caption{Derivatives interval  (shown in grey color) for the function  defined
    by
     formula (\ref{5.6.1})  with $C=0$. }
 \label{Big_paper7}
  \end{center}
\end{figure}

If the chosen numeral system, $\mathcal{S}$, is such that the
point $x=0$   belongs to $[a,b]_{\mathcal{S}}$ then (\ref{5.6.2})
is also true for $x \neq 0$ and we should verify the existence  of
the derivative or derivatives interval for $x=0$. It follows from
(\ref{5.3.4}) and (\ref{5.4.3.0}) that at this point the left
derivative exists and is   calculated as follows
\[
f_1(x)-f_1(x-l) = -x-(-(x-l))= l \cdot (-1),
 \]
 \[
f^{-}(x,l) =   -1, \hspace{1cm} f^{-}(x)= f^{-}(x,0) = -1.
\]
The right derivative is obtained analogously from (\ref{5.3.5})
and (\ref{5.4.3.0})
  \[
f_3(x+r)-f_3(x)  =   x+C+r-(x+C)= r \cdot 1,
\]
\[
f^{+}(x,r) =   1, \hspace{1cm} f^{+}(x)= f^{+}(x,0)=  1.
\]
However, $f(x)$ has continuous formulae at the point $x=0$  only
for $C=0$. Thus, only in  this case  function $f(x)$ defined by
  (\ref{5.6.1}) has the derivatives interval $[-1,1]$ at the
point $x=0$. This situation is illustrated in
Fig.~\ref{Big_paper7}.

By a complete analogy it is possible to show that functions
 \[
 \begin{array}{cc}
\hat{f}(x)  = \left\{ \begin{array}{ll} -2\mbox{\ding{172}}x, &
x \le 0,\\
 3\mbox{\ding{172}}x^2,  & x >0,
 \end{array} \right.   & \hspace{1cm}
 \tilde{f}(x)  = \left\{ \begin{array}{ll} -4\mbox{\ding{172}}^{-1.6}x, &
x \le 0,\\
 5\mbox{\ding{172}}^{-28}x,  & x >0,
 \end{array} \right. \end{array}
\]
have at the point $x=0$ derivatives intervals
$[-2\mbox{\ding{172}},0]$ and
$[-4\mbox{\ding{172}}^{-1.6},5\mbox{\ding{172}}^{-28}]$,
respectively. The first of them is  infinite and the second
infinitesimal.
 \hfill $\Box$
\end{example}

Let us find now the derivative of the function $f(x)=\sin(x)$. As
in standard and non-standard analysis, this is done by appealing
to geometrical arguments.

\begin{lemma}
\label{l1} There exists a function $g(x)$ such that the function
$f(x)=\sin(x)$ can be represented by the following formula:
$\sin(x)=x\cdot g(x),$ where $g(0)=1$.
\end{lemma}

\textit{Proof.} Let us consider a function
$d(x)=\frac{\sin(x)}{x}$ and study it on the right from the point
$x=0$, i.e., in the form $d_3(\xi)$ from (\ref{5.3.2}) for $\xi>0$
(the function $d_1(\xi)$ is investigated by a complete analogy).
It is easy to show from geometrical considerations (see
\cite{Courant}) that
\[
1 < \frac{\xi}{\sin(\xi)} < \frac{1}{\cos(\xi)}, \hspace{1cm} 0 <
\xi < \frac{\pi}{2}.
\]
Inasmuch as $\xi \neq 0$, it follows from these estimates that
\[
\cos(\xi)  < \frac{\sin(\xi)}{\xi} < 1, \hspace{1cm} 0 < \xi <
\frac{\pi}{2},
\]
and, therefore,
\[
\cos(\xi)  <  d_3(\xi)  < 1, \hspace{1cm} 0 < \xi < \frac{\pi}{2},
\]
 Due to Corollary~\ref{c3}, the function $\cos(x)$
has continuous formulae and $\cos(0)=1$. As a consequence, we
obtain $d_3(0)=1$. Since the formula
$d(\xi)=\frac{\sin(\xi)}{\xi}$ has $\xi$ in the denominator, the
only possibility to execute the reduction in the formula
$d_3(\xi)$ leading to the result $d_3(0)=1$ consists of the
existence of the   function $g(\xi)$ such that $d_3(\xi)=\frac{\xi
\cdot g(\xi)}{\xi}=g(\xi)$. The fact that $\sin(x)$ has continuous
formulae concludes the proof. \hfill $\Box$

This lemma illustrates again our methodological point of view
expressed in Postulate~2. The formula $d(x)=\frac{\sin(x)}{x}$ has
zero in the denominator when $x=0$. In contrast, the function
$g(x)$ is defined also for $x=0$ in spite of the fact that they
describe the same mathematical object. The difficulties we have in
$d(x)$ have been introduced by inadequate mathematical instruments
used to describe the object, particularly, by the fact that the
precise formula for $\sin(x)$ is unknown. Analogously, if we
rewrite the constant function $c(x)=1$ in the form $m(x)= x \cdot
\frac{1}{x}$ we obtain the same effect. Introduction of the
designations  (\ref{5.3.2}) allows us to monitor this situation
easily.

\begin{theorem}
\label{tder3} The derivative of the function $f(x)=\sin(x)$ is
$f'(x)=\cos(x)$.
\end{theorem}

\textit{Proof.} Let us calculate the right relative difference,
$f^{+}(x,r)$, from (\ref{5.3.5}) for $f(x)=\sin(x)$ at a point $x$
 as follows
 \[
f_3(x+r)-f_3(x)  =  \sin(x+r)-\sin(x)=
\sin(x)\cos(r)+\cos(x)\sin(r)-\sin(x).
\]
By using the trigonometric identity
\[
\cos(x)=1-2\sin^2(0.5x)
\]
and  Lemma~\ref{l1} we obtain
 \[
f_3(x+r)-f_3(x) =\sin(x)(\cos(r)-1)+\cos(x)\sin(r)=
\]
 \[
\sin(r)\cos(x)-2\sin^2(0.5r)\sin(x)=r\cdot g(r)
\cos(x)-0.5r^2\cdot g(0.5r)^2\sin(x) =
\]
 \[
 r  \big(g(r)  \cos(x)-0.5r \cdot
g(0.5r)^2\sin(x)\big).
\]
Thus, we have   that the right relative difference
  \[
f^{+}(x,r) =  g(r)  \cos(x)-0.5r \cdot g(0.5r)^2\sin(x).
\]
By recalling again Lemma~\ref{l1} we obtain from this formula that
  \[
f^{+}(x) = f^{+}(x,0) =  g(0)  \cos(x) = 1 \cdot \cos(x)=\cos(x).
\]
By a complete analogy (see (\ref{5.3.4})) we obtain that $f^{-}(x)
= \cos(x)$, too. Due to the introduced definition of the
derivative (see (\ref{5.3.4})--(\ref{5.4.3.0}), (\ref{5.4.5})), it
follows from the fact $f^{+}(x) = f^{-}(x)$ that $f'(x)=\cos(x)$.
 \hfill $\Box$

 We conclude this section by an example showing the
usage of derivatives for calculating sums with an infinite number
of items. Recall that due to Postulate~3 and since we have
infinite numbers that can be written explicitly, we cannot define
a function in the form $f(x) = \sum_{i=0}^{\infty} a_i(x)$. It is
necessary to indicate explicitly an infinite number, $n$, of items
in the sum $f(x) = \sum_{i=0}^{n} a_i(x)$ and, obviously, two
different infinite numbers, $n=n_1$ and $n=n_2$, will define two
different functions.
\begin{example}
\label{e24} Let us consider an infinite number $n$ and the
following function
 \beq
f(x) = \sum_{i=0}^{n} x^i = 1 + x + x^2 + \ldots + x^n,
 \label{der7}
 \eeq
where $x\neq 1$ can be finite, infinite or infinitesimal. By
derivating (\ref{der7}) we obtain
 \beq
f'(x) = \sum_{i=1}^{n} i x^{i-1} = 1 + 2x + 3x^2 + \ldots + n
x^{n-1}.
 \label{der8}
 \eeq
However, it follows from (\ref{3.7.2.f}) and (\ref{3.7.2.f.1})
that
 \beq
  f(x) =
\frac{1-x^{n+1}}{1-x}.
 \label{der9}
\eeq
 By derivating (\ref{der9}) we have that
 \beq
  f'(x) =
\frac{1+nx^{n+1}-(n+1)x^{n}}{(1-x)^2}.
 \label{der10}
\eeq
 Thus, we can conclude that for finite and infinite   $n$ and
for   finite, infinite or infinitesimal $x\neq 1$ it follows
 \beq
\frac{1+nx^{n+1}-(n+1)x^{n}}{(1-x)^2} = \sum_{i=1}^{n} i x^{i-1} =
1 + 2x + 3x^2 + \ldots + n x^{n-1},
 \label{der11}
\eeq
 whereas the traditional analysis uses just the following formula
 \[
\frac{1}{(1-x)^2} = \sum_{i=0}^{\infty} i x^{i-1} = 1 + 2x + 3x^2
+ \ldots  ,
\]
that is able to provide a result only if for any infinite $n$
there exists a finite approximation of (\ref{der11}) and it is
used only for finite $x,$ $-1 < x < 1$.
  \hfill $\Box$
\end{example}

\section{Conclusion}
\label{s7}

In this paper, a new applied methodology to Calculus has been
proposed. It has been emphasized that the philosophical triad --
the researcher, the object of investigation, and tools used to
observe the object -- existing in such natural sciences as Physics
and Chemistry exists in Mathematics, too. In natural sciences, the
instrument used to observe the object influences results of
observations. The same happens in Mathematics  studying numbers
and objects that can be constructed by using numbers. Thus,
numeral systems used to express numbers are instruments of
observations used by mathematicians. The usage of powerful numeral
systems gives the possibility of obtaining more precise results in
Mathematics in the same way as usage of a good microscope gives
the possibility to obtain more precise results in Physics.

A brief introduction to a unified theory of continuous and
discrete functions  has been given for functions  that can assume
finite, infinite, and infinitesimal values over finite, infinite,
and infinitesimal domains. This theory allows one to work with
derivatives  and subdifferentials     that can assume finite,
infinite, and infinitesimal values, as well. It has been shown
that the expressed point of view on Calculus allows one to avoid
contrapositions of the previous approaches (finite quantities
versus infinite, standard analysis versus non-standard, continuous
analysis versus discrete, numerical analysis versus pure) and to
create a unique framework for Calculus having a simple and
intuitive structure. It is important to emphasize that the new
approach has its own computational tool -- the Infinity Computer
-- able to execute numerical computations with finite, infinite,
and infinitesimal quantities.

\bibliographystyle{amsplain}
\bibliography{XBib_Derivatives}

\providecommand{\bysame}{\leavevmode\hbox to3em{\hrulefill}\thinspace}
\providecommand{\MR}{\relax\ifhmode\unskip\space\fi MR }
\providecommand{\MRhref}[2]{%
  \href{http://www.ams.org/mathscinet-getitem?mr=#1}{#2}
}
\providecommand{\href}[2]{#2}
\begin{thebibliography}{10}

\bibitem{Benci}
V.~Benci and M.~{Di~Nasso}, \emph{Numerosities of labeled sets: a new way of
  counting}, Advances in Mathematics \textbf{173} (2003), 50--67.

\bibitem{Cantor}
G.~Cantor, \emph{Contributions to the founding of the theory of transfinite
  numbers}, Dover Publications, New York, 1955.

\bibitem{Cauchy}
A.L. Cauchy, \emph{Le calcul infinit\'esimal}, Paris, 1823.

\bibitem{Clarke}
F.H. Clarke, \emph{Optimization and nonsmooth analysis}, Wiley, New York, 1983.

\bibitem{Conway}
J.H. Conway and R.K. Guy, \emph{The book of numbers}, Springer-Verlag, New
  York, 1996.

\bibitem{Courant}
R.~Courant and F.~John, \emph{Introduction to calculus and analysis {I}},
  Springer-Verlag, New York, 1999.

\bibitem{DAlembert}
J.~d'Alembert, \emph{Diff\'erentiel}, Encyclop\'edie, ou dictionnaire
  raisonn\'e des sciences, des arts et des m\'etiers \textbf{4} (1754).

\bibitem{Gordon}
P.~Gordon, \emph{Numerical cognition without words: {E}vidence from
  {A}mazonia}, Science \textbf{306} (2004), no.~15 October, 496--499.

\bibitem{Leibniz}
G.W. Leibniz and J.M. Child, \emph{The early mathematical manuscripts of
  {L}eibniz}, Dover Publications, New York, 2005.

\bibitem{Newton}
I.~Newton, \emph{Method of fluxions}, 1671.

\bibitem{fractals}
H.-O. Peitgen, H.~J\"{u}rgens, and D.~Saupe, \emph{Chaos and fractals},
  Springer-Verlag, New York, 1992.

\bibitem{Pica}
P.~Pica, C.~Lemer, V.~Izard, and S.~Dehaene, \emph{Exact and approximate
  arithmetic in an amazonian indigene group}, Science \textbf{306} (2004),
  no.~15 October, 499--503.

\bibitem{Robinson}
A.~Robinson, \emph{Non-standard analysis}, Princeton Univ. Press, Princeton,
  1996.

\bibitem{Sergeyev}
{Ya.D.} Sergeyev, \emph{Arithmetic of infinity}, Edizioni Orizzonti
  Meridionali, CS, 2003.

\bibitem{Sergeyev_patent}
Ya.D. Sergeyev, \emph{Computer system for storing infinite, infinitesimal, and
  fi\-ni\-te quan\-ti\-ties and executing arithmetical operations with them},
  patent application 08.03.04, 2004.

\bibitem{www}
\bysame, \emph{http://www.theinfinitycomputer.com}, 2004.

\bibitem{Poland}
\bysame, \emph{Mathematical foundations of the {I}nfinity {C}omputer}, Annales
  {UMCS} {I}nformatica {AI} \textbf{4} (2006), 20--33.

\bibitem{chaos}
\bysame, \emph{Blinking fractals and their quantitative analysis using infinite
  and infinitesimal numbers}, Chaos, Solitons $\&$ Fractals \textbf{33(1)}
  (2007), 50--75.

\bibitem{spirals}
\bysame, \emph{Measuring fractals by infinite and infinitesimal numbers},
  Mathematical Methods, Physical Methods $\&$ Simulation Science and Technology
  \textbf{1(1)} (2008), 217--237.

\bibitem{biology}
\bysame, \emph{Modelling season changes in the infinite processes of growth of
  biological systems}, Transactions on Applied Mathematics and Nonlinear Models
  (2008), (to appear).

\bibitem{informatica}
\bysame, \emph{A new applied approach for executing computations with infinite
  and infinitesimal quantities}, Informatica \textbf{19(4)} (2008), 567--596.

\bibitem{Korea}
\bysame, \emph{Numerical computations and mathematical modelling with infinite
  and infinitesimal numbers}, Journal of Applied Mathematics and Computing
  \textbf{29} (2009), 177--195.

\bibitem{Wallis}
J.~Wallis, \emph{Arithmetica infinitorum}, 1656.

\end{thebibliography}
\end{document}